\documentclass[12pt,a4paper]{article}

\usepackage{amsmath,amssymb}

\def\cM{{\cal M}}
\def\oM{{\overline{\cal M}}}
\def\bM{{\partial\cM}}
\def\oMbd{{\cal D}}
\def\oC{{\overline{\cal C}}}
\def\qed{{\hfill $\diamondsuit$}}

\def\CP{{{\mathbb C}{\rm P}}}

\def\Aut{{\rm Aut}}
\def\Z{{\mathbb Z}}

\def\C{{\mathbb C}}
\def\Q{{\mathbb Q}}
\def\d{{\partial}}
\def\cT{{\cal T}}
\def\ident{
\begin{picture}(18,10)
\put(3,-1){$\rightarrow$}
\put(3.5,2.5){$\sim$}
\end{picture}
}
\def\ie{{\it i.e.}}
\def\pt{{\rm pt}}

\def\ev{{\rm ev}}

\def\whs{\widehat{s}}
\def\whu{\widehat{u}}
\def\whv{\widehat{v}}

\usepackage{theorem}

\newtheorem{theorem}{Theorem}
\newtheorem{proposition}{Proposition}[section]
\newtheorem{corollary}[proposition]{Corollary}
\newtheorem{lemma}[proposition]{Lemma}

{\theorembodyfont{\rmfamily}
\newtheorem{definition}[proposition]{Definition}
\newtheorem{example}[proposition]{Example}
\newtheorem{remark}[proposition]{Remark}

}

\usepackage{graphicx}

\title{Tautological relations and the $r$-spin Witten conjecture}

\author{Carel Faber\thanks{Institutionen f{\"o}r Matematik, Kungliga
Tekniska H{\"o}gskolan, S-10044 Stockholm, Sweden. 
E-mail: faber@kth.se.   
},
Sergey Shadrin\thanks{Korteweg\,-\,de Vries Instituut voor
Wiskunde, Universiteit van Amsterdam, Plantage Muidergracht 24,
1018 TV Amsterdam, Netherlands. Department of Mathematics, Institute of
System Research, Nakhimovskii prospekt 36-1, Moscow 117218, Russia.
E-mail: s.shadrin@uva.nl, 
shadrin@mccme.ru\,.
},
Dimitri Zvonkine\thanks{
Institut math{\'e}matique de Jussieu,
Universit{\'e} Paris~VI, 175, rue du Chevaleret,
75013 Paris, France. E-mail: zvonkine@math.jussieu.fr. \newline \newline
{\em 2000 Mathematics Subject Classification:} 
14H10, 14N35, 53D45, 53D50. \newline
{\em Key words:} Quantization of Frobenius manifolds, Gromov--Witten potential,
moduli of curves, $r$-spin structures, Witten's conjecture.
}}

\date{}

\begin{document}

\maketitle

\begin{abstract}
In~\cite{Lee1,Lee2}, Y.-P.~Lee introduced a notion of
{\em universal relation} for formal Gromov--Witten potentials.
Universal relations are connected to tautological
relations in the cohomology ring of the moduli space
$\oM_{g,n}$ of stable curves. Y.-P.~Lee conjectured that
the two sets of relations coincide and proved the inclusion
(tautological relations) $\subset$ (universal relations)
modulo certain results announced by C.~Teleman.
He also proposed an algorithm that, conjecturally,
computes all universal/tautological relations.

Here we give a geometric interpretation of Y.-P.~Lee's
algorithm. This leads to a much simpler proof of the
fact that every tautological relation gives rise to
a universal relation. We also show that Y.-P.~Lee's
algorithm computes the tautological relations correctly
if and only if the Gorenstein conjecture on the tautological cohomology 
ring of $\oM_{g,n}$ is true. 
These results are first steps in the task of establishing
an equivalence between formal and geometric Gromov--Witten
theories.

In particular,
it implies that in any semi-simple Gromov--Witten theory where arbitrary
correlators can be expressed in genus~0 correlators using
only tautological relations, the formal and the geometric
Gromov--Witten potentials coincide.

As the most important application we show that our results
suffice to deduce the statement of a 1991 Witten conjecture
on $r$-spin structures from the results obtained by
Givental for the corresponding formal Gromov--Witten potential.

The conjecture in question states that certain intersection numbers
on the moduli space of $r$-spin structures can be arranged
into a power series that satisfies the $r$-KdV (or $r$th higher
Gelfand-Dikii) hierarchy of partial differential equations.
\end{abstract}

\tableofcontents

\section{Introduction}

\subsection{Tautological relations and universal relations}

Tautological cohomology classes of $\oM_{g,n}$ are, roughly, intersections
of boundary components with $\psi$-classes and $\kappa$-classes.
They are described by so-called {\em dual graphs}
(Definition~\ref{Def:dual}). Sometimes it happens that
a nontrivial linear combination $L$ of dual graphs determines
a zero cohomology class $[L]=0$. Such linear combinations are called
{\em tautological relations}. Until now there is no
known way to generate all of them by an algorithm.

Every tautological relation gives rise to a large family of
relations in every geometric Gromov--Witten theory. Indeed,
let $X_{g,n+n',D} = \oM_{g,n+n'}(X,D)$ be the space of
stable maps from genus~$g$ curves with~$n+n'$ marked points
to a target compact K{\"a}hler manifold~$X$, of degree 
$D \in H_2(X,\Z)$. Denote by $[X_{g,n+n',D}]$ its virtual
fundamental class. (Virtual fundamental classes are defined in
\cite{BehFan,LiTian,Siebert}, but we will actually use only
certain factorization properties of these classes.) 
Let $p:X_{g,n+n',D} \rightarrow \oM_{g,n}$ be the 
forgetful map. Further, let~$L$ be a tautological relation on~$\oM_{g,n}$,
and $\beta \in H^*(X_{g,n+n',D})$ a cohomology class of the
form $\prod_{i=1}^{n+n'} \psi_i^{d_i} \ev_i^*(\mu_i)$, where
$\mu_i$ are cohomology classes of~$X$ and $\ev_i$ the evaluation maps.
Then the integral
\begin{equation} \label{Eq:integral}
\int\limits_{[X_{g,n+n',D}]} \!\!\!\!\! p^*[L]\; \beta
\end{equation}
can be expressed as a polynomial in Gromov--Witten invariants of~$X$
in a standard way. On the other hand, we know in advance
that this integral vanishes because $[L]=0$.

In~\cite{Givental}, A.~Givental introduced a notion of
{\em formal Gromov--Witten descendant potentials}. It is conjectured that
formal Gromov--Witten potentials contain all possible
geometric Gromov--Witten potentials and more. The relation
between formal Gromov--Witten theory and geometry
is not established, although C.~Teleman announced some important
results on this subject\footnote{Teleman's preprint~\cite{Teleman}
with a classification of semi-simple cohomological field theories
appeared a year after the first version of this paper was completed.}. 
Nonetheless, since there exists
a universal expression in Gromov--Witten 
correlators that represents integral~(\ref{Eq:integral})
for geometric Gromov--Witten potentials, it is natural to
conjecture that this expression vanishes not only
in the geometric case, but for {\em all formal Gromov--Witten
potentials} as well.

This led Y.-P.~Lee to introduce the notion
of a {\em universal relation} in~\cite{Lee1,Lee2}. Roughly speaking, it
is a linear combination of dual graphs such that the
induced relations for Gromov--Witten correlators hold for all formal
Gromov--Witten potentials, see Definition~\ref{Def:univrel}.
Y.-P.~Lee conjectured that the tautological relations and
the universal relations were exactly the same and
suggested an algorithm to compute them. To put things
clearly, let Taut denote the set of tautological relations,
Univ the set of universal relations, and Alg the set of
relations computed by Y.-P.~Lee's algorithm. 

The inclusion $\mbox{Alg} \subset \mbox{Univ}$ is a reformulation of
Theorem~\ref{Thm:Lee}. It was proved by Y.-P.~Lee 
in~\cite{Lee2} and we re-explain the proof here in
Section~5 with some additional details.

Y.-P.~Lee also stated as a theorem the inclusion
$\mbox{Taut} \subset \mbox{Univ}$, although his proof
is quite complicated and involves a reference to
C.~Teleman's results.

In Section~\ref{Ssec:Gorenstein} we prove the following theorem.

\begin{theorem} \label{Thm:alggoren}
A linear combination of dual graphs~$L$ belongs to $\mbox{\rm Alg}$
if and only if the intersection numbers of the class $[L]$
with all tautological classes of complementary dimension
vanish.
\end{theorem}

It follows immediately that 
$\mbox{Taut} \subset \mbox{Alg} \subset \mbox{Univ}$.
Moreover, we have $\mbox{Taut} = \mbox{Alg}$ if and
only if the so-called {\em Gorenstein conjecture},
stating that the Poincar{\'e} pairing is nondegenerate on
the tautological cohomology ring of $\oM_{g,n}$, is true.

We have no new information about the conjectured equalities 
$\mbox{Alg} \stackrel{?}{=} \mbox{Univ}$ and 
$\mbox{Taut}\stackrel{?}{=} \mbox{Univ}$.

\bigskip

Summarizing, we see that if a linear combination~$L$ of
dual graphs represents a zero cohomology class, then
all formal Gromov--Witten potentials ``know about it''
in the sense that the expressions in correlators that represent
intersection numbers with the pull-back of~$[L]$ vanish.

This result constitutes a rather powerful 
tool for establishing equivalences between formal and geometric
Gromov--Witten potentials:

{\em In any semi-simple Gromov--Witten theory where any correlator
can be expressed in terms of genus~0 correlators using only
tautological relations, the formal Gromov--Witten potential
coincides with the geometric Gromov--Witten potential.}

A more precise formulation is given in Theorem~\ref{Thm:GeomForm}.

As the main application, we show that
this result suffices to deduce Witten's conjecture on the
space of $r$-spin structures from the results already established
by Givental for formal Gromov--Witten potentials.

\subsection{Witten's conjecture}
\label{Ssec:Witten}

In 1991 E.~Witten formulated a conjecture on the intersection
theory of the spaces of $r$-spin structures of Riemann 
surfaces~\cite{Witten}. An $r$-spin structure on a smooth
curve~$C$ with $n \geq 1$ marked points $x_1, \dots, x_n$
is a line bundle~$\cal T$ together with an identification
$$
\cT^{\otimes r} \ident K(-\sum a_i x_i),
$$
where $K$ is the cotangent line bundle and
the integers $a_i \in \{0, \dots, r-1\}$ are chosen
in such a way that $2g-2-\sum a_i$ is divisible by~$r$.
The space of $r$-spin structures has a natural compactification
$\oM^{1/r}_{g; a_1, \dots, a_n}$ with a forgetful map
$p:\oM^{1/r}_{g; a_1, \dots, a_n} \rightarrow \oM_{g,n}$
(see~\cite{Jarvis,AbrJar,Chiodo1}). 
Using the line bundle~$\cT$ and the map~$p$ one can define
on $\oM_{g,n}$ a cohomology class $c_W(a_1, \dots, a_n)$
of pure dimension, now called {\em Witten's class}. 
{\em Assuming} that $H^0(C,\cT) = 0$ for every stable 
curve~$C$, we find that $V = H^1(C, \cT)$ is a vector 
bundle over $\oM^{1/r}_{g; a_1, \dots, a_n}$. Then 
Witten's class is defined by
$$
c_W(a_1, \dots, a_n) = \frac1{r^g} p_*(c_{\rm top} (V^\vee)).
$$
In other words: take the dual vector bundle of~$V$, take its
Euler (or top Chern) class, take its push forward to
$\oM_{g,n}$, and divide by $r^g$. In the general case,
when $H^0(C,\cT)$ does not vanish identically, Witten's
class has a much more intricate definition that we do not
recall here.

Witten's conjecture
states that the intersection numbers of Witten's class
with powers of $\psi$-classes can be arranged into a
generating series satisfying the $r$-KdV (or $r$th
higher Gelfand-Dikii) hierarchy.

\begin{theorem} \label{Thm:Witten}
Witten's $r$-spin conjecture is true.
\end{theorem}

Below we give a summary of our proof of this conjecture.
The ideas we use are close to those used by Y.-P.~Lee
to prove the conjecture in low genus in~\cite{Lee3,Lee4}.

\begin{enumerate}
\item
In his initial paper 
E.~Witten~\cite{Witten} proved that the conjecture was
true in genus~0 provided the class~$c_W$ satisfied certain postulated
properties. The space of $r$-spin structures and the
class~$c_W$ were later constructed precisely and shown to 
possess the expected properties by a joint effort of several 
people~\cite{Jarvis,JaKiVa1,Mochizuki,PolVai,Polishchuk}.

\item
A.~Givental~\cite{Givental} 
constructed a transitive group action on all semi-simple
formal Gromov--Witten theories. He found a specific 
group element that takes
the Gromov--Witten potential of a point to the string solution
of the $r$-KdV hierarchy. (In particular, the genus zero part of this
solution coincides with the genus zero part of the generating
function for Witten's correlators.)

\item
Y.-P.~Lee~\cite{Lee1, Lee2} found an algorithm that allows one to compute
universal relations on Gromov--Witten potentials in 
Givental's theory, \ie, relations satisfied by all formal 
Gromov--Witten descendant potentials.

\item 
In this paper we give a geometric interpretation of Y.-P.~Lee's
algorithm and deduce a simple proof of Y.-P.~Lee's claim
that any tautological relation in
$H^*(\oM_{g,n})$ gives rise to a universal relation
in Givental's theory.

\item
The first author together with R.~Pandharipande~\cite{FabPan}
studied the so-called double ramification classes in $H^*(\oM_{g,n})$
and proved that they were tautological. Therefore integrals of
Witten's class and powers of $\psi$-classes on these cycles
correspond to some polynomials in correlators
of the Gromov--Witten potential in Givental's setting.

\item
Finally, the last two authors~\cite{ShaZvo} showed that
tautological relations on double ramification classes were
sufficient to express Witten's correlator in terms of genus 
zero correlators.

\end{enumerate}

The logic of our proof is the following. Item~6 gives an
algorithm that expresses any given correlator in terms of genus~0
correlators. The steps of this algorithm can be performed in
Givental's theory by Item~4. Thus the generating series for
the correlators coincides with the Gromov--Witten potential
given by Givental. But the latter is known to give a 
solution of the $r$-KdV hierarchy by Item~2.

Useful references on Witten's $r$-spin conjecture also include
\cite{Kontsevich,AbrJar,CaCaCo,JaKiVa2,JaKiVa3,JaKiVa4,
Chiodo1,Chiodo2,Shadrin2}.

\section{Formal Gromov--Witten 
potentials and tautological classes}
\label{Sec:Givental}

In this section we give a short overview of the aspects
of Givental's theory of formal Gromov--Witten potentials that we
will need.

\subsection{Formal Gromov--Witten descendant potentials}

\begin{definition}
The {\em genus~$g$ Gromov--Witten descendant potential of a point} 
is the formal power series
$$
F_g^\pt(t_0, t_1, \dots) = 
\sum_{n \geq 0} \, \sum_{d_1, \dots,  d_n} \!\!
\left<
\tau_{d_1} \dots \tau_{d_n}
\right>
\frac{t_{d_1} \dots t_{d_n}}{n!} ,
$$
where
$$
\left<
\tau_{d_1} \dots \tau_{d_n}
\right> =
\int\limits_{\oM_{g,n}} \psi_1^{d_1} \dots \psi_n^{d_n},
$$
$g$ being determined by the condition 
$\dim \oM_{g,n} = 3g-3+n = \sum d_i$.
The {\em total Gromov--Witten descendant potential of a point}
is $F^\pt = \sum F^\pt_g \hbar^{g-1}$ and its exponential 
$Z^\pt = \exp F^\pt$ is called the {\em Gromov--Witten partition
function} of a point.
\end{definition}

A formal Gromov--Witten potential is defined to model
certain properties of $F^\pt$ and those of
Gromov--Witten potentials of more general target K{\"a}hler 
manifolds~$X$. We restrict our considerations to the even part of the 
cohomology of $X$. In our description we explain in brackets 
the geometric aspects that motivate the axiomatic definitions.

Let $V$ be a complex vector space [the space $H^{\rm even}(X,\C)$]
in which we choose a basis~$A$. The space~$V$
possesses a distinguished element~$1$ 
[cohomology class~1], which we usually assume to be the first
vector of the basis. $V$ also carries a nondegenerate symmetric bilinear 
form~$\eta$ [Poincar{\'e} pairing]. The coefficients of $\eta$
in the basis will be denoted by $\eta_{\mu\nu}$ and the
coefficients of the inverse matrix by $\eta^{\mu\nu}$.

\begin{definition}
Let $M$ be a neighborhood of the origin in~$V$.
Let $F_0$ be a power series in variables $t_d^\mu$,
$d =1,2,3, \dots$, $\mu \in A$, whose coefficients are 
analytic functions on~$M$ in variables $t_0^\mu$. The coefficients
of $F_0$ are denoted by
$$
F_0 = 
\sum_{n \geq 0} \, 
\sum_{
\stackrel{\scriptstyle d_1, \dots,  d_n \ge 0}
{\mu_1, \dots, \mu_n \in A}
} 
\!\!
\left<
\tau_{d_1, \mu_1} \dots \tau_{d_n,\mu_n}
\right>
\frac{t_{d_1}^{\mu_1} \dots t_{d_n}^{\mu_n}}{n!} .
$$
$F_0$ is called a {\em formal genus~$0$ Gromov--Witten potential} if
it satisfies the string equation, the dilaton equation, and
the topological recursion relation
(see, for instance,~\cite{Givental}, \cite{Givental3}
or~\cite{Givental4}). The open set~$M$ is called a
{\em Frobenius manifold}\/\footnote{Sometimes the definition
of a Frobenius manifold also includes an Euler field
and $F_0$ is required to satisfy certain homogeneity
conditions with respect to this field. In other
sources Frobenius manifolds with an Euler field are
called {\em conformal}.}.
\end{definition}

\begin{remark} \label{Rem:Novikov}
If $X$ is a K{\"a}hler manifold and $E \subset H_2(X,\Z)$ 
its semi-group of effective $2$-cycles, one usually considers 
Gromov--Witten potentials with coefficients not in $\C$, but
in the {\em Novikov ring} of power series of the form
$$
\sum_{D \in E} c_D Q^D, \qquad c_D \in \C.
$$
We will mostly work with Gromov--Witten potentials over~$\C$,
since we do not need Novikov rings
for our main application, namely the Witten conjecture.
However, we indicate in remarks
the modifications that must be made in the general case. 
A detailed introduction to formal Gromov--Witten potentials,
including a discussion of Novikov rings, can be found 
in~\cite{LeePan}.
\end{remark}

Given a formal genus~0 Gromov--Witten potential~$F_0(t_d^\mu)$, 
let $f_0(t_0^\mu)$ be the series obtained from~$F_0$
by setting $t_d^\mu = 0$ for $d \geq 1$. On every tangent
space to~$M$ one defines an algebra via the structural constants
$$
C^{\mu_3}_{\mu_1,\mu_2} (t_0^\mu)
= \sum_{\nu} 
\frac{\d^3 f_0}{\d t_0^{\mu_1} \d t_0^{\mu_2} \d t_0^\nu} 
\; \eta^{\nu, \mu_3}
$$
depending on $t_0^\mu \in M$. Together with the bilinear form~$\eta$,
one gets the structure of a {\em Frobenius algebra} in every
tangent space to~$M$.

\begin{definition} \label{Def:semisimple}
$F_0$ is called {\em semi-simple} if the algebra structure
is semi-simple for generic $t_0^\mu \in M$. The {\em rank} of~$F_0$
is the dimension of~$V$.
\end{definition}

\begin{remark}
A (genus~$0$) Gromov--Witten potential defined over a Novikov
ring~$R$ is called semi-simple if the algebra structure
at a generic point is semi-simple over the
{\em algebraic closure of the field of fractions} of~$R$. 
\end{remark}

\begin{example}
Let $F$ be the string solution of the 3-KdV hierarchy
(the Gromov--Witten potential that appears in
Witten's conjecture for $r=3$, see Section~\ref{Ssec:WiGrWi}). 
Denote $x= t^0_0$, $y=t^1_0$. Then
$$
f_0(x,y) = \frac{x^2y}2 + \frac{y^4}{72}.
$$
Let $F$ be the Gromov--Witten potential of $\CP^1$.
Denote $x= t^1_0$, $y=t^\omega_0$, where $1\in H^0(\CP^1)$ and 
$\omega \in H^2(\CP^1)$ form the natural basis of $H^*(\CP^1)$. Then
$$
f_0(x,y) = \frac{x^2y}2 + Q e^y.
$$
\end{example}

\bigskip

In~\cite{Givental3}, Givental constructs an action of
the so-called {\em twisted loop group} on
Gromov--Witten potentials of rank~$k$.
More precisely, the group itself is ill-defined, but
there is a well-defined action of its upper triangular
and lower triangular parts. 

This action is almost transitive on semi-simple
potentials of rank~$k$. Denote by
\begin{equation} \label{Eq:rescaled}
F^{\pt, \alpha} = \sum \hbar^{g-1} \, \alpha^{2-2g-\sum d_i}
\left< \tau_{d_1} \dots \tau_{d_n} \right> 
\frac{t_{d_1} \dots t_{d_n}}{n!} 
\end{equation}
a rescaled Gromov--Witten potential of the point.

Then, whenever $F_0$ is semi-simple and of rank~$k$, 
there exists an element~$S$ of the lower-triangular
group and an element $R$ of the upper-triangular group
such that $SR$ takes
$F^{\pt, \alpha_1}_0 \oplus \dots \oplus F^{\pt, \alpha_k}_0$ 
to $F_0$, where $\alpha_1, \dots, \alpha_k$ are appropriately
chosen constants.

A quantized version of the same group (described
in~\cite{Givental}), that we will call
{\em Givental's group}, acts on Gromov--Witten 
descendant partition functions. 
The quantization $\widehat{S}(\hbar)\widehat{R}(\hbar) $ of~$SR$ 
acts on power series in~$\hbar$ and $t_d^\mu$. It takes 
$Z^{\pt, \alpha_1} \times \dots \times Z^{\pt, \alpha_k}$ 
to some power series $Z = \exp \sum \hbar^{g-1} F_g$. 
More precisely, $Z$ is in general a power series in variables $t_d^\mu$
for $d \geq 1$, whose coefficients are functions on~$M$
in variables $t_0^\mu$. These functions are analytic on~$M$
outside of the discriminant, \ie, the subvariety where
the tangent Frobenius algebra is not semi-simple.

In practice it often happens that for mysterious reasons
the coefficient functions of~$Z$ turn out to be regular 
on the discriminant of~$M$. In particular, this is the case
for the $A_r$ singularity, which is the main example of
interest for us, and for Gromov--Witten potentials of
K\"ahler manifolds. In the sequel we will only consider
the cases where $Z$ is regular at the origin and thus can
be decomposed into a power series. 
However, the examples of~\cite{DubZha}, Section~6, 
show that singularities on the discriminant do appear in 
many cases.

Taking the lowest degree term in~$\hbar$ (that is, the coefficient
of $\hbar^{-1}$) in the logarithm of~$Z$, we recover the
action of~$SR$.

In Section~5 we will give a precise description of the
upper and lower triangular groups and their actions.

\begin{definition} \label{Def:genusext}
The series $F= \ln Z$ is called a {\em genus 
expansion} of the genus zero potential~$F_0$ or
a {\em formal Gromov--Witten descendant potential.} 
\end{definition}

In terminology coming from physics, the Taylor coefficients
$
 \left<
\tau_{d_1, \mu_1} \dots \tau_{d_n,\mu_n}
\right>_{g,D}
$
of the Gromov--Witten potential~$F$ are called 
{\em correlators}, while the elements 
of the basis~$A$ are called {\em primary fields}.

\medskip

The genus expansion of a formal genus~$0$ potential is not unique,
because the action of the twisted loop group is not free.
Indeed, the upper-triangular subgroup of 
the {\em rank~1} twisted loop group
acts trivially on $F^\pt_0$. The direct product of $k$ copies
of the rank~1 upper-triangular subgroup forms a subgroup
of the rank~$k$ upper triangular subgroup. This subgroup acts trivially on 
$F^{\pt, \alpha_1}_0 \oplus \dots \oplus F^{\pt, \alpha_k}_0$. 
But the quantizations of the elements of this subgroup 
do not act trivially on 
$Z^{\pt, \alpha_1} \times \dots \times Z^{\pt, \alpha_k}$,
and thus we obtain several different genus expansions.

However, the ambiguity of the genus expansion can be fixed
using an additional property of {\em homogeneity}.

To illustrate this property, let us first 
assume that $Z =\exp \sum \hbar^{g-1} F_g$ 
is the descendant Gromov--Witten partition function
of a target K{\"a}hler manifold~$X$. 
Choose a homogeneous basis $A$ of $H^*(X)$. 
To each variable $t_d^{\mu}$ we assign its {\em weight} 
$w(t_d^{\mu}) = d+\deg(\mu)-1$,
where $\deg(\mu)$ is the algebraic degree of~$\mu \in A$. Further, 
introduce a weight function on the Novikov ring:
$w(Q^D) =  -\left< D, c_1(TX) \right>$ for an effective divisor~$D$.
Finally, denote by $\dim X$ the dimension of~$X$.

The expected dimension of~$X_{g,n,D}$ is 
$n + (1-g) (\dim X -3) + \left< D, c_1(TX) \right>$.
Therefore the correlator
$
 \left<
\tau_{d_1, \mu_1} \dots \tau_{d_n,\mu_n}
\right>_{g,D}
$
vanishes unless
$$
\sum d_i + \sum \deg(\mu_i) = n + (1-g) (\dim X -3) + 
\left< D, c_1(TX) \right>
$$
$$
\Longleftrightarrow
\sum w(t_{d_i}^{\mu_i}) + w(Q^D) = (1-g) (\dim X -3).
$$
In other words, $F_g$ is quasihomogeneous of total weight
$(1-g)(\dim X -3)$.

This property is formalized in the following definition.

\begin{definition} \label{Def:homog}
Introduce a weight map $w:A \rightarrow \Q$ such that
$w(1) = -1$, let $w(t_d^\mu) = d+w(\mu)$. Also introduce
a weight function on the Novikov ring.
Let $\dim \in \Q$ be a constant. 

A formal Gromov--Witten potential $F=\sum \hbar^{g-1} F_g$ 
is called  {\em homogeneous}
with respect to $w$ and $\dim$ if every $F_g$ is quasihomogeneous
of total weight~$(1-g)(\dim-3)$.
\end{definition}

A homogeneous formal genus~$0$ Gromov--Witten potential has
a {\em unique} homogeneous genus expansion~(\cite{Givental}
Proposition~6.7~c and Remark~6.9~a).


Givental conjectured that if~$F_0$ is the geometric genus~0 potential of
a target space~$X$, then the total geometric Gromov--Witten potential
of~$X$ coincides with the homogeneous genus expansion of its genus~$0$
part.

\subsection{Tautological classes and dual graphs}

\subsubsection{The $\kappa$-classes}

\begin{definition}
The cohomology class $\kappa_{k_1, \dots, k_m}$ on
$\oM_{g,n}$ is defined by
$$
\kappa_{k_1, \dots, k_m}=
\pi_*(\psi_{n+1}^{k_1+1} \dots \psi_{n+m}^{k_m+1}),
$$
where
$$
\pi:\oM_{g,n+m} \rightarrow \oM_{g,n}
$$
is the forgetful map.
\end{definition}

This definition is compatible with the usual definition
of $\kappa$-classes $\kappa_k$ (for $m=1$). The classes
$\kappa_{k_1, \dots, k_m}$ and the monomials
$\kappa_{k_1} \dots \kappa_{k_m}$ form two bases of the
same vector space and the matrix of basis change is triangular.
Indeed, we have
$$
\kappa_{k_1, \dots, k_m}= \sum_{\sigma \in S_m} \, 
\prod_{c = \mbox{\scriptsize cycle of } \sigma} 
\!\!\!\!\!\!\! \kappa_{k(c)}, 
\quad \mbox{where} \quad
k(c) = \sum_{i \in c } k_i.
$$
For instance, 
\begin{eqnarray*}
\kappa_{k_1,k_2} &=& \kappa_{k_1} \kappa_{k_2} + \kappa_{k_1+k_2},\\
\kappa_{k_1,k_2,k_3} &=& \kappa_{k_1} \kappa_{k_2} \kappa_{k_3} + 
\kappa_{k_1+k_2} \kappa_{k_3} + 
\kappa_{k_1+k_3} \kappa_{k_2} + 
\kappa_{k_2+k_3} \kappa_{k_1} + 
2 \kappa_{k_1+k_2+k_3}.
\end{eqnarray*}
We prefer to work with the classes $\kappa_{k_1, \dots, k_m}$
because they are easier to express in terms of Gromov--Witten
correlators. 

Let $p: \oM_{g,n+1} \rightarrow \oM_{g,n}$ be the
forgetful map.
Let $r : \oM_{g-1,n+2} \rightarrow \oM_{g,n}$ and
$q: \oM_{g_1, n_1+1} \times \oM_{g_2, n_2+1} 
\rightarrow \oM_{g,n}$ be the usual ``gluing'' mappings to the
boundary components of $\oM_{g,n}$. 

\begin{lemma} \label{Lem:kappa}
We have
\begin{eqnarray*}
p^*(\kappa_{k_1, \dots, k_m}) &=& 
\kappa_{k_1, \dots, k_m}-\sum_{i=1}^m \psi_{n+1}^{k_i} \, 
\kappa_{k_1, \dots, \widehat{k_i}, \dots, k_m},\\
q^*(\kappa_{k_1, \dots, k_m}) &=& \!\!\!\!\!
\sum_{I \sqcup J = \{1, \dots, m \} } \!\!\!\!\!
\kappa_{k_I} \times \kappa_{k_J},\\
r^*(\kappa_{k_1, \dots, k_m}) &=& \kappa_{k_1, \dots, k_m},
\end{eqnarray*}
where $k_I = \{k_i \}_{i \in I}$, $k_J = \{k_i \}_{i \in J}$,
and $\widehat{k_i}$ means that the index is omitted.
\end{lemma}

\paragraph{Proof.} Only the first equality is nontrivial.
Consider the forgetful map
$$
{\widetilde p} : \oM_{g,n+1+m} \rightarrow \oM_{g,n+m}.
$$
To avoid confusion in indices, suppose the $n+m$ marked points
are numbered from 1 to~$n+m$, while the forgotten point is
labeled with~$\alpha$. In $H^2(\oM_{g,n+m+1})$ we have
${\widetilde p}\,^*(\psi_{n+i}) = \psi_{n+i} - D_{n+i, \alpha}$,
where $D_{n+i, \alpha}$ is the divisor of curves on which the points
$n+i$ and $\alpha$ lie on a separate sphere with no other marked
points. From the relations
$$
\psi_{n+i} D_{n+i, \alpha} =0 \quad \mbox{for} \quad 1 \leq i \leq m,
\qquad 
D_{n+i, \alpha} D_{n+j, \alpha} =0 \quad \mbox{for} \quad i \not= j,
$$
we obtain that
$$
{\widetilde p}\,^*
\left(\prod_{i=1}^m \psi_{n+i}^{k_i+1} \right)
= 
\prod_{i=1}^m \psi_{n+i}^{k_i+1} 
+
\sum_{i=1}^m \psi_{n+1}^{k_1+1} \dots (-D_{n+i,\alpha})^{k_i+1}
\dots \psi_{n+m}^{k_m+1}.
$$
Taking the push-forward of this class to $\oM_{g,n+1}$
we obtain the right-hand side of the first equality of the lemma.

The second equality comes from the fact that
each of the $m$ points forgotten by the map
$\pi:\oM_{g,n+m} \rightarrow \oM_{g,n}$ can find
itself on either of the two components of the boundary curves. 
The third equality follows directly from the definition. \qed

\subsubsection{The dual graphs}

Consider a stable curve~$C$ of genus~$g$ with~$n$ marked
points in a one-to-one correspondence with a set of markings~$S$.
The topological type of~$C$ can be described
by a graph~$G$ obtained by replacing every irreducible component
of the curve by a vertex and every node of the curve by an edge.
Every marked point is replaced by a {\em tail}
(an edge that does not lead to any vertex)
retaining the same marking as the marked point. 
Each vertex~$v$ is labeled by an integer
$g_v$: the geometric genus of the corresponding component. The
$g_v$'s and the first Betti number of~$G$ add up to~$g$.

In order to avoid problems with automorphisms, we will
label all the half-edges of~$G$. To $G$ we assign the space 
$$
\oM_G = \prod_v \oM_{g_v,n_v},
$$
where the product goes over the set of vertices of~$G$,
$g_v$ is the genus of the vertex~$v$, and $n_v$ its valency
(the number of half-edges and tails adjacent to it). The space $\oM_G$
comes with a natural map $p:\oM_G \rightarrow \oM_{g,n}$
whose image is the closure of the
set of stable curves homeomorphic to~$C$. Note that
$p_*[\oM_G] = |\Aut(G)| \cdot [p(\oM_G)]$.

We can define a cohomology class on 
$\oM_G$ (and hence on $\oM_{g,n}$ taking the push-forward by~$p$)
by assigning a class $\kappa_{k_1, \dots, k_m}$ to each vertex
of~$G$ and a power $\psi^d$ of the $\psi$-class to each half-edge
and each tail of~$G$.

\begin{definition} \label{Def:dual}
A graph~$G$ with labeled half-edges and tails,
describing the topological type of a stable curve,
with an additional label $\kappa_{k_1, \dots, k_m}(v)$
assigned to each vertex~$v$ and
a nonnegative integer~$d_e$ assigned to each half-edge
and each tail~$e$ is called a {\em stable dual graph}.
The corresponding cohomology class of~$\oM_{g,n}$ is called
the {\em tautological class} assigned to~$G$ and denoted
by $[G]$. The {\em genus} of a dual graph is the genus of
the corresponding stable curves, its
{\em degree} is the algebraic degree of the corresponding
cohomology class, and its {\em dimension} is 
$\dim \oM_{g,n} - \mbox{degree} = 3g-3+n-\mbox{degree}$. 
\end{definition}

This definition allows several modifications. 

First of all, we can consider {\em not necessarily connected
dual graphs}. They represent tautological classes on direct
products of several moduli spaces. 

Second, we can consider dual graphs that describe the 
topology of a {\em semi-stable}
curve that is not necessarily stable. In this case we will
also label each tail of the graph with a primary field $\mu \in A$.
If we work over a Novikov ring, we also assign an 
effective 2-cycle $D \in E$ to the whole
graph. Such a graph will be called a {\em semi-stable dual graph with
primary fields}. It describes a tautological cohomology class in 
the space of stable maps $\oM_{g,n+n'}(X,D)$.

A ``dual graph'' with no other specifications
will mean ``stable connected dual graph''. 

\subsection{Universal relations}
\label{Ssec:univrel}

Let $F= \sum \hbar^{g-1} F_g$ be the geometric Gromov--Witten
potential of some target K{\"a}hler manifold~$X$, $A$
a basis of $H^*(X)$ and $\eta$ the Poincar{\'e} pairing.
Let $L = \sum c_i G_i$ be a linear combination 
of dual graphs representing a class $[L] \in H^*(\oM_{g,n})$. 
As in the introduction, let $X_{g,n+n',D} = \oM_{g,n+n'}(X,D)$, where
$D \in E$ is an effective 2-cycle, 
let $[X_{g,n+n',D}]$ be its virtual fundamental class,
and let $p:X_{g,n+n',D} \rightarrow \oM_{g,n}$ be the forgetful
map. We are going to describe a
way to express the integrals of the form
$$
\sum_{D \in E}\; Q^D \int\limits_{[X_{g,n+n',D}]}
\!\!\!\!\!\!
p^*([L]) \; \prod_{i=1}^{n+n'} \psi_i^{d_i} \ev_i^*(\mu_i)
$$
via $\eta^{\mu\nu}$ and the coefficients of the series
$F_0, \dots, F_g$.

We start with assigning a polynomial in correlators
and coefficients $\eta^{\mu\nu}$ to any stable dual
graph or any semi-stable dual graph with primary fields
that does not contain $\kappa$-classes.

\begin{definition} \label{Def:P_G}
Let $G$ be a stable dual graph or a semi-stable dual graph with
primary fields such that no $\kappa$-classes are assigned to its
vertices. We define the {\em polynomial~$P_G$} by the following
procedure. (i)~Assign a primary field $\mu \in A$ to every
half-edge of~$G$. In the case of a stable dual graph, assign,
moreover, the distinguished primary field $1 \in A$ to the tails.
(ii)~To every vertex~$v$ assign the correlator
$\left<\tau_{d_1,\mu_1} \dots \tau_{d_{n_v},\mu_{n_v}} \right>_{g_v}$,
where $g_v$ is the genus of~$v$, $n_v$ its valency, and
$d_i, \mu_i$ the labels on the half-edges and tails adjacent to~$v$.
(iii)~To every edge assign the coefficient $\eta^{\mu\nu}$, where
$\mu$ and $\nu$ are the primary fields corresponding to its half-edges.
(iv)~Take the product of all the correlators and the coefficients
$\eta^{\mu\nu}$ thus obtained. (v)~Sum over all the ways to
attribute primary fields to the half-edges.
If we work over a Novikov ring, we must also sum over all
the ways to assign effective $2$-cycles $D_v$ to the vertices
in such a way that $\sum D_v = D$.
\end{definition}

Note that every edge of~$G$ introduces a contraction of indices 
via the bilinear form~$\eta$. This comes from the fact that the
class of the diagonal in $X \times X$ equals 
$\sum_{\mu,\nu} \eta^{\mu\nu} \, \mu \times \nu$.

Note also that the definition works perfectly well for not 
necessarily connected dual graphs.

\bigskip

Now we go back to our problem of constructing a pull-back
in $X_{g,n+n',D}$ of the tautological class~$[L]$.
Let $G$ be a stable dual graph participating in the linear
combination~$L$.

\paragraph{Step 1: eliminating the $\kappa$-classes.}
If a vertex $v$ of $G$ is labeled with $\kappa_{k_1, \dots, k_m}$,
erase this label and replace it by $m$ new tails issuing
from~$v$ with labels $\psi^{k_1+1}, \dots, \psi^{k_m+1}$
on them. Thus we obtain a new dual graph~$G_1$. The newly
added tails will be called {\em $\kappa$-tails}.

This rule is justified by the following remark: if
$\pi: \oM_{g,n+m} \rightarrow \oM_{g,n}$ is the forgetful map,
we have
$$
\psi_1^{d_1} \dots \psi_n^{d_n} \kappa_{k_1, \dots, k_m}
= \pi_*(\psi_1^{d_1} \dots \psi_n^{d_n}
\psi_{n+1}^{k_1+1} \dots \psi_{n+m}^{k_m+1})
$$
in $H^*(\oM_{g,n})$.

\begin{example}
Let
\begin{center}
\
\begin{picture}(150,30)
\put(-35,10){$G=$}
\put(0,10){\includegraphics{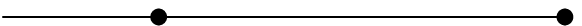}}
\put(32,0){$g=2$}
\put(150,0){$g=1$}
\put(41,20){$\kappa_1$}
\put(15,20){$\psi^0$}
\put(65,20){$\psi^0$}
\put(135,20){$\psi^1$}
\end{picture}
\end{center}
(here $g=3$, $n=1$). Then we have
$$
P_{G_1}=
\sum_{\mu,\nu \in A}
\left<\tau_{0,1} \tau_{2,1} \tau_{0,\nu}
\right>_2
\; \eta^{\nu\mu} \;
\left<\tau_{1,\mu}
\right>_1.
$$
\end{example}

\paragraph{Step 2: recomputing the $\psi$-classes.}
There is a difference
between the $\psi$-classes on $\oM_{g,n}$ and on
$X_{g,n+n',D}$, because of the presence of additional
marked points and because of the appearance of semi-stable
source curves. To take this into account, we modify $G_1$
according to the following rule: replace every half-edge 
and every tail with label $\psi^d$
\begin{center}
\
\begin{picture}(60,15)
\put(0,0){\includegraphics{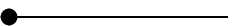}}
\put(30,7){$\psi^{d}$}
\end{picture}
\end{center}
by the linear combination
\begin{center}
\
\begin{picture}(150,30)
\put(0,10){\includegraphics{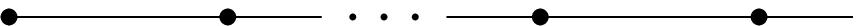}}
\put(-100,10){$\sum\limits_{p=0}^d \,(-1)^p \!\!\!\!
\sum\limits_{d_0 + \dots + d_p = d-p}$}
\put(63,17){$0$}
\put(216,17){$0$}
\put(153,17){$0$}
\put(196,17){$\psi^0$}
\put(132,17){$\psi^0$}
\put(41,17){$\psi^0$}
\put(228,17){$\psi^{d_0}$}
\put(170,17){$\psi^{d_1}$}
\put(80,17){$\psi^{d_{p-1}}$}
\put(15,17){$\psi^{d_p}$}
\end{picture}
\end{center}
Thus we obtain a linear combination of dual graphs~$G_2$. When
we perform Steps~1 and~2 with all graphs of the linear
combination~$L$, we obtain a new linear combination~$L_2$ of semi-stable
dual graphs (without primary fields). It represents the 
class $p^*([L])$ and will therefore be 
denoted by $L_2 = p^*(L)$.

In terms of polynomials $P_G$ assigned to the dual graphs,
replacing $G_1$ by the linear combination of dual graphs~$G_2$
is equivalent to making the following 
substitutions in~$P_{G_1}$:
\begin{eqnarray*}
\tau_{1,\mu} &\mapsto& 
\tau_{1,\mu} - 
\sum_{\mu_1,\nu_1} 
\left<\tau_{0,\mu} \tau_{0,\nu_1} \right>_0
\; \eta^{\nu_1\mu_1} \; \tau_{0,\mu_1},\\
\tau_{2,\mu} &\mapsto& 
\tau_{2,\mu} - 
\sum_{\mu_1,\nu_1} 
\left<\tau_{1,\mu} \tau_{0,\nu_1} \right>_0
\; \eta^{\nu_1\mu_1} \; \tau_{0,\mu_1}
- \sum_{\mu_1, \nu_1} 
\left<\tau_{0,\mu} \tau_{0,\nu_1} \right>_0
\; \eta^{\nu_1\mu_1} \; \tau_{1,\mu_1}\\
&&+ \sum_{\mu_1,\nu_1,\mu_2,\nu_2} 
\left<\tau_{0,\mu} \tau_{0,\nu_1} \right>_0 \; \eta^{\nu_1\mu_1} \;
\left<\tau_{0,\mu_1} \tau_{0,\nu_2} \right>_0
\; \eta^{\nu_2\mu_2} \; \tau_{0,\mu_2},
\end{eqnarray*}
and so on. In general, every insertion of $\tau_{d, \mu}$
must be replaced by
\begin{equation} \label{Eq:ancestors}
\tau_{d,\mu} + 
\sum_{p=1}^d \,(-1)^p \!\!\!\!\!\!\!\!\!\!
\sum_{
\stackrel{\scriptstyle d_0 + \dots + d_p = d-p}
{\mu_1, \nu_1, \dots, \mu_p, \nu_p}
}
\!\!\!\!\!\!\!\!\!
\left<\tau_{d_0, \mu} \tau_{0,\nu_1} \right>_0
\, \eta^{\nu_1 \mu_1} \, 
\left<\tau_{d_1, \mu_1} \tau_{0,\nu_2} \right>_0
\, \eta^{\nu_2 \mu_2} \dots \, \eta^{\nu_p \mu_p} \;\;
\tau_{d_p, \mu_p}.
\end{equation}
This formula is used in the following way: the symbol $\tau_{d,\mu}$
was part of some correlator $\left< x \right>$; now we put 
$\tau_{d_p,\mu_p}$ in its place, while the other factors
of the formula become factors in front of~$\left< x \right>$.


It is important to note that the procedure of Step~2 for
expressing ancestor $\psi$-classes in terms of descendant
$\psi$-classes is universal. This means that the expression remains
valid for any number of additional marked points~$n'$ and
for any target manifold~$X$. 


\paragraph{Step 3: multiplying by 
$\prod_{i=1}^{n+n'} \psi_i^{d_i} \ev_i^*(\mu_i)$.}
Let $G_2$ be one of the dual graphs involved
in the linear combination~$L_2$. It has two kinds of
tails: the $\kappa$-tails and the ordinary tails numbered
from~1 to~$n$. Add $n'$ more tails numbered from~$n+1$
to~$n+n'$ by attaching them to the vertices of $G_2$
in all possible ways. 

Now, the $\kappa$-tails already
bear labels $\psi^k$. We also label them with the
distinguished primary field~$1$.

The tails from~1 to~$n$ also bear labels $\psi^{d'_i}$,
$1 \leq i \leq n$. We replace $\psi^{d'_i}$ on the
$i$th tail by $\psi^{d'_i +d_i}$ and also label the
$i$th tail with the primary field~$\mu_i$.

The tails from $n+1$ to $n+n'$ have no labels. We label them with
$\psi^{d_i}$ and $\mu_i$, $n+1 \leq i \leq n+n'$.

Thus we obtain a linear combination~$L_3$ of semi-stable
dual graphs with primary fields. This linear combination
will be denoted by 
$L_3 = p^*(L) \prod_{i=1}^{n+n'} \psi_i^{d_i} \ev_i^*(\mu_i)$.

In terms of polynomials $P_G$, replacing $P_{L_2}$ by $P_{L_3}$
is equivalent to (i)~replacing the symbols $\tau_{d_i',1}$
corresponding to the tails from 1 to~$n$ by symbols
$\tau_{d_i'+d_i, \mu_i}$ and (ii)~inserting new symbols
$\tau_{d_i,\mu_i}$ for $n+1 \leq i \leq n+n'$ in the
existing correlators in all possible ways.

The corresponding polynomial in correlators will be denoted by
$$
P_{L_3} = \biggl<
p^*(L) \prod_{i=1}^{n+n'} \psi_i^{d_i} \ev_i^*(\mu_i)
\biggr>_g.
$$
If~$F$ is the geometric Gromov--Witten potential of
some target space~$X$, we have
$$
\biggl<
p^*(L) \prod_{i=1}^{n+n'} \psi_i^{d_i} \ev_i^*(\mu_i)
\biggr>_g
=
\sum_{D \in H_2(X)}\; \int\limits_{[X_{g,n+n',D}]}
\!\!\!\!
p^*([L]) \prod_{i=1}^{n+n'} \psi_i^{d_i} \ev_i^*(\mu_i).
$$
However it will be important for us that the expression
$
\biggl<
p^*(L) \prod_{i=1}^{n+n'} \psi_i^{d_i} \ev_i^*(\mu_i)
\biggr>_g
$ 
makes sense for any power series $F = \sum \hbar^{g-1} F_g$ 
in variables $t_d^{\mu}$,
in particular for formal Gromov--Witten potentials. Note that
if two different linear combinations $L$ and $L'$
of dual graphs represent the same cohomology class,
it is, for the time being, not at all obvious that 
$
\biggl<
p^*(L) \prod_{i=1}^{n+n'} \psi_i^{d_i} \ev_i^*(\mu_i)
\biggr>_g
$ 
and
$
\biggl<
p^*(L') \prod_{i=1}^{n+n'} \psi_i^{d_i} \ev_i^*(\mu_i)
\biggr>_g
$ 
coincide for every formal Gromov--Witten potential.

\begin{definition} \label{Def:vector}
Let $L$ be a linear combination of dual graphs.
Let $F = \sum \hbar^{g-1}F_g$ be a power series in variables
$t_d^{\mu}$. The infinite vector of the values of 
$
\biggl<
p^*(L) \prod_{i=1}^{n+n'} \psi_i^{d_i} \ev_i^*(\mu_i)
\biggr>_g
$
for all $n' \geq 0$ and for all $d_i, \mu_i$,
is called the {\em induced vector} of~$L$ and is
denoted by~$F_L$.
\end{definition}

\begin{definition} \label{Def:univrel}
A linear combination~$L$ of dual graphs is called a 
{\em universal relation} if the vector $F_L$ is equal to $0$ for any 
semi-simple formal Gromov--Witten descendant potential~$F$,
as defined by Givental.
\end{definition}

\begin{proposition}
The vanishing of $F_L$ can be expressed as an
infinite family of partial differential equations
on $F_0$, $F_1$, \dots.
\end{proposition}

This is a standard fact so we do not prove it here but
instead illustrate it with an example.

\begin{example}
The tautological relation 
\begin{center}
\
\begin{picture}(180,30)
\put(-3,-9){\includegraphics{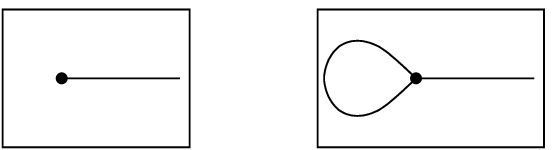}}
\put(0,0){$g=1$}
\put(40,16){$\psi$}
\put(111,0){$g=0$}
\put(59,9.3){$\displaystyle -\frac{1}{12}$}
\put(165, 9.3){$=0$}
\end{picture}
\end{center}
gives rise to the following family of partial differential equations
indexed by $d \geq 0 $ and $\rho \in A$:
$$
\frac{\d F_1}{\d t_{d+1}^\rho} - \sum_{\mu, \nu}
\frac{\d F_1}{\d t_0^\mu} \eta^{\mu \nu} 
\frac{\d^2 F_0}{\d t_0^\nu d t_d^\rho}
- \frac1{24} \sum_{\mu, \nu}
 \frac{\d^3 F_0}{\d t_0^\mu \d t_0^\nu d t_d^\rho}
\eta^{\mu \nu} =0.
$$
\end{example}

\bigskip

The notion of universal relation is naturally generalized to 
not necessarily connected stable dual graphs, but this should
be done carefully. The simplest way to define a universal relation
in this case is the following. 

Consider any graph as a product of its connected components.
Then a linear combination~$L$ of not necessarily connected
stable dual graphs is called a {\em universal relation}
if it can be represented as a sum of several products
such that every term of the sum contains a universal relation
for connected graphs as a factor. 

In other words, if we consider the product 
$\oM_{g_1,n_1} \times \dots \times \oM_{g_k,n_k}$,
a universal relation in one of the factors times any class
in the product of the remaining factors is a universal relation,
and a linear combination of universal relations is still a universal
relation.

\begin{remark} \label{Rem:UnivDisconnected}
The method we used to define the
induced vector $F_L$ works for a linear combination~$L$
of not necessarily connected stable dual graphs. However,
in this case a more natural notion is that of an
{\em extended induced vector} $\widehat{F}_L$. Its definition
is similar to that of $F_L$, with the difference that
when we add new marked points to the curve, we are allowed
to prescribe the connected component to which every point
should go. $L$ is a universal relation if and only if
$\widehat{F}_L$ vanishes. We do not know whether the
conditions $\widehat{F}_L = 0$ and $F_L=0$ are equivalent.
\end{remark}

%
%
%
%
%

Now we can give a precise formulation of the theorem announced 
in the introduction.

Let $F = \sum \hbar^{g-1} F_g$ 
be a power series in $\hbar$ and $t_d^{\mu}$,
where $F_0$ is a formal genus~$0$ Gromov--Witten potential.
We introduce five properties of~$F$.

\begin{enumerate}
\item
{\bf Homogeneity.} $F$ is homogeneous (in the sense of 
Definition~\ref{Def:homog}) for some $w$ and $\dim$.

\item
{\bf Geometricity.} 
$F_L$ vanishes for every tautological relation~$L$.

\item 
{\bf Semi-simplicity.} $F_0$ is semi-simple at some point~$t$.

\item
{\bf Reducibility to genus~$0$.} Every correlator of~$F$ can be
expressed in terms of genus~$0$ correlators using only
properties~1 and~2.

\item
{\bf Analyticity.} The homogeneous genus expansion of~$F_0$
is regular at the origin.
\end{enumerate}

All geometric Gromov--Witten potentials of target K{\"a}hler
manifolds satisfy conditions~1 and~2. Conditions~3, 4, and~5,
on the other hand, must be checked in every particular case.

In Section~\ref{Sec:Witten} we will see an example of a potential
that satisfies all five
conditions without being the
geometric potential of a target manifold.

\begin{theorem} \label{Thm:GeomForm}
A power series~$F$ satisfying conditions~{\rm (1-5)}
coincides with the homogeneous genus expansion of~$F_0$.
\end{theorem}

\paragraph{Proof.} Let $\widehat F$ be the 
homogeneous genus expansion of~$F_0$. The expressions
$F_L$ and $\widehat F_L$ vanish for all tautological 
relations~$L$ (the vanishing of $\widehat F_L$ 
follows from Theorems~\ref{Thm:alggoren} and~\ref{Thm:Lee}). 
Similarly, the nonhomogeneous correlators vanish both for~$F$ 
and for~$\widehat F$. According to condition~4, these vanishing
conditions are enough to express every correlator 
in terms of genus~$0$ correlators. But the genus~0 correlators
coincide, since both are given by~$F_0$.
Thus $F = \widehat F$. \qed

\bigskip

In~\cite{Lee2}, Y.-P.~Lee studied the action of Givental's group
on the induced vectors~$F_L$ and constructed an algoritm that
computes certain, conjecturally all, universal relations.
We are now going to state his results.

\subsection{Y.-P. Lee's algorithm}

\subsubsection{The operators $\tau_k$}
\label{Sssec:tau}

We are going to define linear operators $\tau_k$ acting on
the space of linear combinations of dual graphs. 
Here~$k$ is an arbitrary positive integer.

If~$G$ is a connected dual graph whose tails are labeled by a set~$S$,
then $\tau_k(G)$ is a linear combination of not necessarily
connected dual graphs with labeling set~$S \cup \{\alpha,\beta \}$. The
graphs of this linear combination are obtained from~$G$ by the
following operations.

\begin{enumerate}
\item
Cut an edge of~$G$ into two tails. Change their labels to $\alpha$ and
$\beta$ in both possible ways. If the labels on the half-edges
were $\psi^a$ (on~$\alpha$) and $\psi^b$ (on~$\beta$)
before the cutting, we now label them
first with $\psi^{a+k}$ and $\psi^b$ and then with $\psi^a$ and
$\psi^{b+k}$. The first stable graph thus
obtained is taken with coefficient~1 while the second
is taken with coefficient~$(-1)^{k-1}$.

\begin{center}
\ 
\begin{picture}(360,40)
\put(0,15){\includegraphics{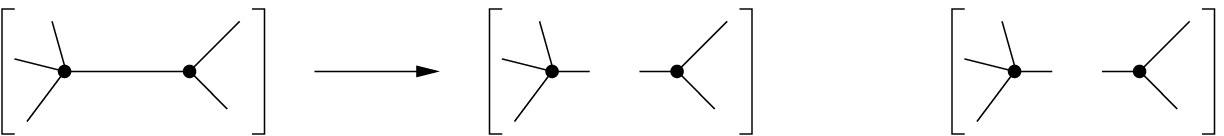}}
\put(22,40){$\psi^a$}
\put(45,40){$\psi^b$}
\put(160,25){$\alpha$}
\put(185,23){$\beta$}
\put(295,25){$\alpha$}
\put(320,23){$\beta$}
\put(160,40){$\psi^{a+k}$}
\put(187,40){$\psi^b$}
\put(294,40){$\psi^a$}
\put(316,40){$\psi^{b+k}$}
\put(220,30){$+ \; (-1)^{k-1}$}
\put(142,0){$+$ two more terms.}
\end{picture}
\end{center}

\item
Split a vertex~$v$ of~$G$ in two, and add a new tail on each of them,
one marked by~$\alpha$ and the other one by~$\beta$. If the genus
of~$v$ was~$g$, assign to the new vertices genera $g_1$ 
and $g_2$ such that $g_1 + g_2 = g$ in all possible ways. Distribute
the edges that were going out of~$v$ between the two vertices
in all possible ways. If~$v$ carried the label $\kappa_{k_1, \dots, k_m}$,
split the set $\{ k_1, \dots, k_n \}$ in two disjoint subsets
$I$ and $J$ in all possible ways and assign to the new vertices
the labels $\kappa_I$ and $\kappa_J$. Label tail~$\alpha$
with $\psi^i$ and tail~$\beta$ with $\psi^j$, in all possible ways
with the condition $i+j = k-1$. The dual graph thus obtained is
taken with coefficient $(-1)^{j+1}$. Keep only stable
graphs and sum over all the possibilities
described above.
\begin{center}
\ 
\begin{picture}(280,135)
\put(-32,0){\includegraphics{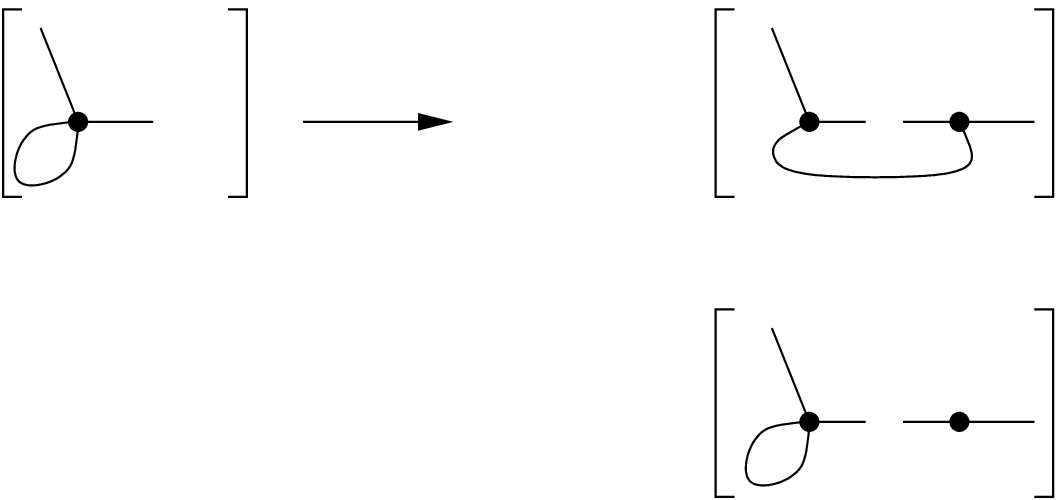}}
\put(-6,117){$g$}
\put(-6,100){$\kappa_{k_1, \dots, k_m}$}
\put(108,106){\LARGE  $\sum$}
\put(133,108){$(-1)^{j+1}$}
\put(185,110){$g_1$}
\put(245,116){$g_2$}
\put(193,99){$\kappa_I$}
\put(250,103){$\kappa_J$}
\put(210,103){$\alpha$}
\put(228,100){$\beta$}
\put(210,115){$\psi^i$}
\put(228,115){$\psi^j$}
\put(97,23){$+$}
\put(108,21){\LARGE  $\sum$}
\put(133,23){$(-1)^{j+1}$}
\put(187,26){$g_1$}
\put(245,31){$g_2$}
\put(189,14){$\kappa_I$}
\put(245,14){$\kappa_J$}
\put(210,16){$\alpha$}
\put(228,13){$\beta$}
\put(210,28){$\psi^i$}
\put(228,28){$\psi^j$}
\put(280,23){$+ \dots$}
\end{picture}
\end{center}

\item
Choose a vertex of~$G$, decrease its genus by~1 and add two new tails on it,
one marked by~$\alpha$ and the other one by~$\beta$. Label tail~$\alpha$
with $\psi^i$ and tail~$\beta$ with $\psi^j$, in all possible ways
with the condition $i+j = k-1$. The dual graph thus obtained is
taken with coefficient $(-1)^{j+1}$. Sum over all possible $i$ and~$j$.
\begin{center}
\ 
\begin{picture}(220,55)
\put(-33,0){\includegraphics{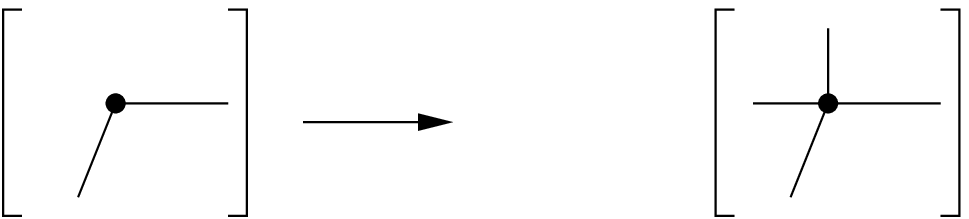}}
\put(-1,23){$g$}
\put(106,24){\LARGE $\sum$}
\put(130,26){$(-1)^{j+1}$}
\put(205,23){$g-1$}
\put(180,26){$\alpha$}
\put(208,50){$\beta$}
\put(180,37){$\psi^i$}
\put(193,50){$\psi^j$}
\end{picture}
\end{center}

\end{enumerate}

The operations~$\tau_k$ are extended to linear combinations
of dual graphs by linearity.

\subsubsection{The algorithm}
\label{Ssec:alg}

\begin{theorem}[Y.-P.~Lee] \label{Thm:Lee}
A linear combination $L$ of dual graphs is
a universal relation whenever
{\rm (i)}~$F^\pt_L = 0$ and
{\rm (ii)}~$\tau_k(L)$ is a universal relation for all~$k\ge1$.
\end{theorem}

This theorem is actually an algorithm for computing universal 
relations for Gromov--Witten potentials. Indeed, the vector
$F^\pt_L$ has an infinite number of entries,
but they can be expressed via a finite number of entries
using the string and dilaton equations. Therefore it is enough
to check Condition~(i) for a finite number of entries.
Now, the dimension~(see Definition~\ref{Def:dual})
of $\tau_k(L)$ is smaller than
that of~$L$, therefore we can proceed by
induction on the dimension of the relations.

Y.-P. Lee conjectures that this algorithm finds
{\em all} universal relations on formal Gromov--Witten
potentials and that these universal relations arise
from tautological relations in the tautological ring
of~$\oM_{g,n}$. However, neither of these claims is proved.

Y.-P.~Lee has also observed that in many cases Condition~(i) is
unnecessary and conjectured that checking Condition~(ii) is
enough for $\dim L\geq 1$. 
However, it follows from our geometric interpretation that this
conjecture is wrong. The first example where it fails is in
$\oM_{2,1}$ in dimension~$1$.

Y.-P.~Lee's proof of Theorem~\ref{Thm:Lee} is summarized
in Section~\ref{Sec:Lee}.

\section{Geometrical interpretation}
\label{Sec:geom}

In this section we give a geometric interpretation of
the operators $\tau_k$ as intersections with the boundary.
(This interpretation was discovered independently by
Y.-P.~Lee and R.~Pandharipande - private communication.)
Using it, we prove that a linear combination~$L$ of dual graphs
is obtained as a relation in Y.-P.~Lee's algorithm if and
only if the intersection of~$[L]$ with all tautological classes
of complementary dimension vanishes.

\subsection{Operators $\tau_i$ and boundary classes}
\label{Ssec:tau}

In the universal curve~$\oC_{g,n}$, consider the
codimension~2 subvariety~$\Delta$ of nodes in the singular
fibers of~$\oC_{g,n}$. Each point of~$\Delta$ is a node of a stable
curve and we will label by~$\alpha$ and~$\beta$
the two marked points of its normalization identified at this node. 
This can be done in two ways, hence
we obtain a double covering of~$\Delta$ that we will
call~$\oMbd$. The space~$\oMbd=\oMbd_{g,n}$ has one connected component
isomorphic to~$\oM_{g-1,n+2}$ (unless~$g=0$) and $(g+1)2^n-2(n+1)$
other connected components, each of which is
isomorphic to~$\oM_{g_1,n_1+1} \times \oM_{g_2,n_2+1}$
for suitable~$g_i$ and~$n_i$
with~$g_1+g_2 = g$ and~$n_1 + n_2 = n$. It comes with a natural map
$p:\oMbd \rightarrow \oM_{g,n}$, whose image is the boundary
$\bM_{g,n}=\oM_{g,n}\setminus\cM_{g,n}$.
Note that all {\em tautological} classes supported on~$\bM_{g,n}$
are defined as push-forwards of classes on~$\oMbd$ under~$p$.

On~$\oMbd$ we define the following cohomology classes:
\begin{eqnarray*}
\rho_1 &=& 1\\
\rho_2 &=& \psi_\alpha - \psi_\beta\\
\rho_3 &=& \psi_\alpha^2 - \psi_\alpha \psi_\beta + \psi_\beta^2\\
\rho_4 &=& \psi_\alpha^3 - \psi_\alpha^2 \psi_\beta +
\psi_\alpha \psi_\beta^2 - \psi_\beta^3
\end{eqnarray*}
and so on.

\begin{proposition} \label{Prop:tauboundary}
Let~$G$ be a dual graph of genus~$g$ with $n$ tails. Then we have
$$
[\tau_k(G)] = - \rho_k \, p^*[G]
$$
in the cohomology ring of~$\oMbd$.
\end{proposition}

\begin{corollary}
Let~$L$ be a linear combination of dual graphs
and suppose that the class~$[\tau_1(L)]$ vanishes
(respectively, has zero intersection with all tautological classes
of complementary dimension).
Then the class~$[\tau_k(L)]$ vanishes 
(respectively, has zero intersection with all tautological classes
of complementary dimension) for all~$k$.
\end{corollary}

\paragraph{Proof.}
The class~$[\tau_k(L)]$ is obtained
from~$[\tau_1(L)]$ by a multiplication by~$\rho_k$.~\qed

\bigskip

This corollary confirms Y.-P.~Lee's experimental observation that 
requiring $\tau_k(L)$ to be a universal relation for all $k \geq 1$
is equivalent to requiring just $\tau_1(L)$ to be a universal relation
(\cite{Lee1}, Section~2.2, Remark~(iii)).

\paragraph{Proof of Proposition~\ref{Prop:tauboundary}.}
The main idea of the proof is very simple. On every boundary
component of $\oM_{g,n}$ we can define the classes $\psi_\alpha$ and 
$\psi_\beta$ corresponding to the marked points $\alpha$ and
$\beta$ identified at the node.
It is well-known that the first Chern class of the
normal line bundle to the boundary
component in $\oM_{g,n}$ equals $-(\psi_\alpha+\psi_\beta)$.

Now, when we intersect a tautological class with a boundary 
component of $\oM_{g,n}$ two cases can occur: either the class is 
entirely contained in the component, or it intersects it transversally.
In the first case we must multiply our class by
the first Chern class of the normal line bundle
$-(\psi_\alpha+\psi_\beta)$ and then by $-\rho_k$. 
Their product is equal to $\psi_\alpha^k + (-1)^{k-1} \psi_\beta^k$. 
In the second case, we must add a new node, that either separates 
a component of the curve in two or is a nonseparating node. In both 
cases we multiply the class thus obtained by 
$-\rho_k = \sum_{i+j = k-1} (-1)^{j+1} \psi_\alpha^i \psi_\beta^j$.
These three possibilities correspond
to the three operations in the definition of~$\tau_k$. 

Now we present the proof with all necessary details.

To begin with, consider the case where~$G$ is a dual graph 
of genus~$g$ with $n$ tails
with trivial labels assigned to all vertices, tails, 
and half-edges. Thus~$[G]$ is the cohomology class
$p_{G*}[\oM_G]$ on~$\oM_{g,n}$.

The computation of $p^*[G]$ fits in the framework of~\cite{GraPan}, A.4.
The contributing graphs are of two kinds: either they have one more edge
than~$G$ or they have the same number of edges as~$G$.

Consider a contributing graph with one extra edge. Contracting the
extra edge, we obtain the graph~$G$ and the contracted edge determines
a unique vertex~$v$ of~$G$; let~$g_v$ be the geometric
genus of the corresponding component and let~$n_v$ be the
valence of~$G$ at~$v$ including both half-edges and tails.
Unless~$g_v=0$, there is exactly one
contributing graph corresponding to~$v$ whose extra edge is
a self-edge at a vertex with geometric genus~$g_v-1$. Moreover, there are
exactly~$(g_v+1)2^{n_v}-2(n_v+1)$ other contributing graphs 
corresponding to~$v$
whose extra edge connects two distinct vertices. In other words, the
contributing graphs whose extra edge contracts to~$v$ are in one-to-one
correspondence with the connected components of the double
covering~$\oMbd_{g_v,n_v}$. Let~$\Gamma$ be a contributing graph with extra
edge~$e$. Let~$h_{\alpha}$ and~$h_{\beta}$ be the labeled half-edges
constituting~$e$ and let~$\alpha$ and~$\beta$ be the corresponding
tails arising upon cutting~$e$.
If cutting~$e$ disconnects the graph, put
$$\oM_e:=\oM_{g_{\alpha},n_{\alpha}}\times\oM_{g_{\beta},n_{\beta}}$$
where~$g_{\alpha}$ resp.~$g_{\beta}$ and~$n_{\alpha}$ resp.~$n_{\beta}$
are the genus and the number of tails of the connected component
containing~$\alpha$ resp.~$\beta$. If~$e$ is non-disconnecting, put
$$\oM_e:=\oM_{g-1,n+2}.$$
Further, let
$$p_{\Gamma,e}:\oM_{\Gamma}\to\oM_e\subset\oMbd_{g,n}$$
be the map obtained via cutting~$e$ and contracting all other edges.
Then by~\cite{GraPan}, A.4 the contribution of~$(\Gamma,e)$ to~$p^*[G]$
equals $$p_{\Gamma,e*}[\oM_{\Gamma}].$$

Next, consider a contributing graph with the same number of edges as~$G$.
The graph may then be identified with~$G$ and exactly one edge~$e$ of~$G$
has been selected: it is identified by the excess line bundle.
Define~$\oM_e$ as above and let
$$p_{G,e}:\oM_G\to\oM_e\subset\oMbd_{g,n}$$
be the map obtained via cutting~$e$ and contracting all other edges.
Then by~\cite{GraPan}, A.4 the contribution of~$(G,e)$ to~$p^*[G]$ equals
$$p_{G,e*}(-\psi_{\alpha}-\psi_{\beta})=
(-\psi_{\alpha}-\psi_{\beta})p_{G,e*}[\oM_G]$$
(with a slight abuse of notation).

We find the following formula for~$p^*[G]$:
$$p^*[G]=\sum_{(\Gamma,e):\,\Gamma/e=G} p_{\Gamma,e*}[\Gamma]
+\sum_{e=(h_{\alpha},h_{\beta})\in E(G)}
(-\psi_{\alpha}-\psi_{\beta})p_{G,e*}[G].
$$
If~$v_0(G)$ denotes the number of vertices of~$G$ of geometric genus~$0$,
the number of summands in the first sum equals
$$\sum_{v\in V(G)}\big((g_v+1)2^{n_v}-2n_v-2\big)+|V(G)|-v_0(G),$$
while the second sum runs over the half-edges of~$G$.

Having dealt with the case of a dual graph with trivial labels, we
consider the case where~$G$ is a decorated dual graph: each vertex~$v$
of~$G$ has been assigned a
class~$\kappa_{K_v}$, where~$K_v$ is a collection of~$m_v$
nonnegative integers,
and each half-edge or tail~$h$ has been assigned a power~$\psi_h^{d_h}$
of its~$\psi$-class~$\psi_h$.

Each half-edge or tail on~$G$ determines a unique such on~$\Gamma$. Taking
this into account, the formula above for~$p^*[G]$ continues to hold
in the case of decorated half-edges or tails.

Let~$v$ be a vertex of~$G$. In case the extra edge~$e$ on~$\Gamma$ doesn't
contract to~$v$, a unique vertex~$w$ of~$\Gamma$ corresponds to~$v$ and it
is decorated with the corresponding class~$\kappa_{K_v}$.
We proceed analogously in case~$e$ is a self-edge at a vertex~$w$
(with genus~$g_w=g_v-1$) contracting to~$v$. Finally, if~$e$ contracts
to~$v$ and connects two distinct vertices, the
decoration~$\kappa_{K_v}$ has to be divided up
in all~$2^{m_v}$ possible ways over the two vertices
(cf.~Lemma~\ref{Lem:kappa}). The formula above for~$p^*[G]$ continues to hold
for an arbitrary decorated dual graph; note that the number
of summands in the first sum equals
$$
\sum_{v\in
V(G)}2^{m_v}\big((g_v+1)2^{n_v}-2n_v-2\big)+|V(G)|-v_0(G).
$$

These considerations prove the statement for~$\tau_1$.
The statement for~$\tau_k$ follows from the trivial
remark that
$$
\rho_k (\psi_\alpha + \psi_\beta) = \psi_\alpha^k +
(-1)^{k-1} \psi_\beta^k.
$$
\qed

\subsection{The algorithm and the Gorenstein conjecture}
\label{Ssec:Gorenstein}

The {Gorenstein conjecture}~\cite{Pandharipande} states that the
tautological ring of $\oM_{g,n}$ has the form of the cohomology
ring of a smooth manifold, that is, it has a unique top degree
class and a nondegenerate Poincar{\'e} duality. 
The conjecture is made for the~$\Q$-subalgebra~$R^*(\oM_{g,n})$ of the
rational Chow ring generated by the tautological classes. In this case
the first claim holds (see~\cite{GraVak,FabPan}), but little
is known about the second claim. In this paper we consider the
Gorenstein conjecture for the tautological {\em cohomology\/}
ring~$RH^*(\oM_{g,n})$, the image in cohomology of~$R^*(\oM_{g,n})$. 
The first claim is then
obvious: the top degree class is the class of a point and it
can be represented by any dual graph of maximal degeneration without
$\psi$- and $\kappa$-classes. 
The second claim is again not known.

\begin{definition}
We say that a tautological cohomology class on $\oM_{g,n}$
is {\em Gorenstein vanishing} if its intersection numbers
with all tautological classes of complementary dimension vanish.
\end{definition}

Here we are going to prove Theorem~\ref{Thm:alggoren}
stating that a linear combination~$L$ of dual graphs
appears as a relation in Y.-P.~Lee's algorithm if and
only if the class~$[L]$ is Gorenstein vanishing. 

\paragraph{Proof of Theorem~\ref{Thm:alggoren}.}
Denote by $p:\oM_{g,n+n'} \rightarrow \oM_{g,n}$ the forgetful map
and by $q:\oMbd \rightarrow \oM_{g,n}$ the natural projection.

Taking into account the form of Y.-P. Lee's algorithm 
(Section~\ref{Ssec:alg}), and the geometric interpretation
of the operators~$\tau_k$, the theorem can be 
reformulated as follows:

{\em A tautological class $\alpha \in H^*(\oM_{g,n})$ 
is Gorenstein vanishing if and only if 

(i)~all the intersection numbers
$$
\int\limits_{\oM_{g,n+n'}} p^*(\alpha) \prod_{i=1}^{n+n'} \psi_i^{d_i},
$$
vanish; and

(ii)~the classes $\rho_k\,  q^*(\alpha)$ are Gorenstein vanishing
for all~$k$.
}

Since $\rho_1 = 1$, the second condition can, of course, be
replaced by the condition $q^*(\alpha)$ is Gorenstein vanishing.

{\bf The only if part.} Suppose $\alpha$ is Gorenstein vanishing.
Then
$$
\int\limits_{\oM_{g,n+n'}} \!\!\!\! p^*(\alpha) \prod_{i=1}^{n+n'}
\psi_i^{d_i}
= 
\int\limits_{\oM_{g,n}} \!\! \alpha \; 
p_* \biggl(\prod_{i=1}^{n+n'} \psi_i^{d_i}\biggr)
=0,
$$
because the class $p_* \left(\prod_{i=1}^{n+n'} \psi_i^{d_i}\right)$
is tautological.

Similarly, if $\beta$ is a tautological class on $\oMbd$,
then
$$
\int\limits_{\oMbd} q^*(\alpha) \; \beta
=
\int\limits_{\oM_{g,n}} \alpha \; q_*(\beta) = 0,
$$
because $q_*(\beta)$ is tautological.

{\bf The if part.} Suppose $\alpha$ satisfies Conditions~(i)
and~(ii), and let $G$ be a dual graph of genus~$g$ with $n$ tails,
of complementary dimension to~$\alpha$.
We wish to prove that $\alpha \cap [G] = 0$.

First consider the case when $G$ has no edges. This means that
$[G]$ is a polynomial in $\psi$- and $\kappa$-classes. A class
like that can be represented as a linear combination of several
classes $p_* \left(\prod_{i=1}^{n+n'} \psi_i^{d_i}\right)$.
But, by the same equality that we used in the ``if'' part,
we have
$$
\int\limits_{\oM_{g,n}} \!\! \alpha \; 
p_* \biggl(\prod_{i=1}^{n+n'} \psi_i^{d_i}\biggr)
\; 
=
\int\limits_{\oM_{g,n+n'}} \!\!\!\! p^*(\alpha) \prod_{i=1}^{n+n'}
\psi_i^{d_i} = 0
$$
by Condition~(i).

Now suppose $G$ has at least one edge. Then the class $[G]$
is the push-forward of some tautological class~$\beta$ on
$\oMbd$, in other words, $[G] = q_*(\beta)$. Thus
$$
\alpha \cap [G]
=
\int\limits_{\oM_{g,n}} \alpha \; q_*(\beta)
=
\int\limits_{\oMbd} q^*(\alpha) \; \beta
 = 0,
$$
because $q^*(\alpha)$ is Gorenstein vanishing by Condition~(ii).
\qed

\section{A proof of Witten's conjecture}
\label{Sec:Witten}

In this section we explain in more detail the plan of the
proof of Witten's conjecture outlined in the introduction.

\subsection{Witten's conjecture and Gromov--Witten theories}
\label{Ssec:WiGrWi}

The generating series studied in Witten's conjecture is
$$
F^{[r]} = \sum_{
\stackrel{\scriptstyle g \geq 0}{n \geq 1}
} \, 
\sum_{
\stackrel{\scriptstyle d_1, \dots, d_n}{a_1, \dots, a_n}
}
\; \int\limits_{\oM_{g,n}} \!\!\!
c_W(a_1, \dots, a_n) \, \psi_1^{d_1} \dots \psi_n^{d_n}
\; \cdot \; 
\frac{t_{d_1}^{a_1} \dots t_{d_n}^{a_n}}{n!},
$$
where $c_W(a_1, \dots, a_n)$, also depending on $r$ and $g$,
is the Witten class briefly introduced in the introduction.

Although this series is not the Gromov--Witten potential of
any target space, it is part of the framework of
Gromov--Witten theories, because its genus~0
part satisfies the string, the dilaton, and the TRR
equations (see~\cite{JaKiVa1}).

The set of primary fields~$A$, its distinguished element, 
and the bilinear form~$\eta$ (that participates in 
the TRR equation) are determined by the following
properties of Witten's class, proved in~\cite{Polishchuk}.

{\bf 1.} If $a_i = r-1$ for some $i$, then $c_W(a_1, \dots, a_n)=0$.
Thus the set $A$ equals $\{ 0, \dots, r-2 \}$.

{\bf 2.} Let $p:\oM_{g,n+1} \rightarrow \oM_{g,n}$ be the
forgetful map. Then 
$$
p^*(c_W(a_1, \dots, a_n)) = c_W(a_1, \dots, a_n ,0).
$$ 
It follows that the distinguished primary field is~0.
(There is an unfortunate clash of notation with the
Givental theory, where the distinguished primary field is usually
denoted by~1, a convention that we followed in our paper
except in the applications to Witten's conjecture.)

{\bf 3.} Let $r : \oM_{g-1,n+2} \rightarrow \oM_{g,n}$ and
$q: \oM_{g_1, n_1+1} \times \oM_{g_2, n_2+1} 
\rightarrow \oM_{g,n}$ be the gluing mappings to the
boundary components of $\oM_{g,n}$. In the second case
we assume, for simplicity, that the $n_1$ marked points
on the component of genus~$g_1$ have numbers $1, \dots, n_1$,
while the $n_2$ marked points
on the component of genus~$g_2$ have numbers $n_1+1, \dots, n$.
Then we have
$$
q^*(c_W(a_1, \dots, a_n)) = \sum_{a'+a''=r-2}
c_W(a_1, \dots, a_{n_1} ,a') \times c_W(a'',a_{n_1+1}, \dots, a_n),
$$
(where at most one term of the sum is actually nonzero, because of
the condition 
$2g_1 - 2 - \sum_{i=1}^{n_1} a_i -a'
 \; \equiv \; 
2g_2 - 2 - \sum_{i=n_1+1}^n a_i -a'' 
\; \equiv \; 0 \; \bmod r$), and
$$
r^*(c_W(a_1, \dots, a_n)) = \sum_{a'+a''=r-2}
c_W(a_1, \dots, a_n ,a',a'').
$$

It follows that the bilinear form is given by
$$
\eta^{ab} = \delta^{a+b,r-2}.
$$

The string, the dilaton, and the TRR equations allow
us to express all genus zero correlators of the
series~$F^{[r]}$ using only the correlators for $g=0$, $n=3,4$, 
which were computed in Witten's original paper.

In~\cite{Givental2}, Givental found a specific element
of Givental's group that sends the series
$$
Z^\pt \times \dots \times Z^\pt \quad
(r-1 \, \mbox{ factors})
$$
to the series $Z^{(r)} = \exp \sum F_g^{(r)} \hbar^{g-1}$ 
such that $Z^{(r)}|_{\hbar =1}$ is the $\tau$-function 
of the string solution of the $r$-KdV hierarchy, while
$F^{(r)} = \sum F_g^{(r)}$ is the string solution itself.
The genus~0 part of $F^{(r)}$ coincides
with the genus~0 part of~$F^{[r]}$. Indeed, both satisfy the
string, the dilaton, and the TRR equations and it is easy to
check that they have the same coefficients for $n=3,4$.

Thus what remains to be proved is that the geometric series~$F^{[r]}$
and the formal Gromov--Witten potential~$F^{(r)}$ coincide
also in higher genus. To do that, we will use 
Theorem~\ref{Thm:GeomForm}, so let us check the conditions~(1-5)
involved in its formulation.

{\bf Homogeneity.} The correlator
$$
\left< 
\tau_{d_1, a_1} \dots \tau_{d_n,a_n}
\right>_g
$$
vanishes unless
$$
\sum d_i + \sum \frac{a_i}r = 
n + \left(\frac{r-2}r -3 \right) (1-g).
$$
Indeed, from the definitions of Section~\ref{Ssec:Witten},
we get $c_1(\cT) = (2g-2 - \sum a_i)/r$, hence the rank of~$V$
and the degree of~$c_W$ equal
$$
\deg c_W = \frac{(r-2)(g-1) + \sum a_i}r.
$$

Thus the series $F^{[r]}$ satisfies the homogeneity property
for the weight function $w(i) = \frac{i}r-1$, and the ``dimension''
$\dim = \frac{r-2}{r}$. (There is no need for a weight function
on a Novikov ring, because we are working over~$\C$) 

{\bf Geometricity.} It is obvious that the Gromov--Witten
potential $F^{[r]}$ respects all tautological relations
because of its geometric origin.

{\bf Semi-simplicity.} There exists a bi-polynomial
isomorphism 
$$
\C^{r-2} \rightarrow \C^{r-2} : 
(t_{0,1}, \dots, t_{0,r-2}) 
\mapsto
(s_1, \dots, s_{r-2})
$$
such that the algebra described in
Definition~\ref{Def:semisimple} is naturally identified
with the algebra 
$$
\C[X]/(X^{r-1} + s_1 X^{r-3} + \dots + s_{r-2}).
$$
Thus it is semi-simple whenever the polynomial has $r-1$ distinct
roots. 

For instance, for $r=5$, the algebra of
Definition~\ref{Def:semisimple} is isomorphic to
$\C[X]/P'$, where 
$$
P(X) = X^5 - t_{0,3} X^3- t_{0,2} X^2  + (t_{0,3}^2/5 - t_{0,1})X
$$
and $P'$ is its derivative.

For more details, see~\cite{Givental2}.

{\bf Reducibility to genus~$0$.} This property will be
established in Theorem~\ref{Thm:rspin}.

{\bf Analyticity.} In~\cite{Givental2} the genus
expansion of $F^{(r)}_0$ is given in the form of a
power series at the origin.

\subsection{Admissible covers and double ramification cycles}

The spaces of admissible covers and the double ramification
cycles were first introduced by Ionel~\cite{Ionel} and
proved very useful in the study of moduli spaces. Let us
briefly recall their definitions.

Consider a map $\varphi$ from a smooth curve $C$ with $n$
marked points to the sphere $S = \CP^1$. On $S$ we mark 
all branch points of~$\varphi$ and the images of the $n$
marked points of~$C$. On~$C$ we then mark all the preimages of the
points that are marked on~$S$.
Now choose several disjoint simple loops on $S$, that do
not pass through the marked points. Suppose that if we contract
these loops we obtain a stable genus~0 curve $S'$. Now contract
also all the preimages of the loops in~$C$ to obtain a
nodal curve $C'$ that turns out to be automatically stable. 
We have obtained a map $\varphi'$ from a stable curve~$C'$ 
of genus~$g$ to a stable curve~$S'$ of genus~0. It has the
same degree over every component of $S'$. Moreover, at each
node of $C'$, the projection $\varphi'$ has the same local
multiplicity on both components meeting at the node.

\begin{definition} \label{Def:AdmCov}
A map from a stable curve of genus~$g$ to a stable curve of
genus~0 topologically equivalent to a map described above is called
an {\em admissible covering}. 
\end{definition}

We will be particularly interested in the space of admissible
coverings with multiple ramifications over only~2 points labeled
with $0$ and~$\infty$, the other ramification points being simple.

\begin{definition} \label{Def:DR-Space}
Consider the space of admissible coverings of some given genus~$g$
with prescribed ramification types 
over two points labeled $0$ and $\infty$, and with
simple ramifications elsewhere. The normalization of this space
is called a {\em double ramification space} or a {\em DR-space}. 
\end{definition}

\begin{definition} \label{Def:DR-cycle} Let $k_1, \dots, k_{n+p}$
be a list of integers such that $\sum k_i =0$ and $k_i \not=0$
for $n+1 \leq i \leq n+p$.
Consider the set of smooth curves $(C, x_1, \dots, x_n) \in \cM_{g,n}$
such that there exist $p$ more marked points $x_{n+1}, \dots, x_{n+p}$
and a meromorphic function on $C$ with no zeroes or poles outside of
$x_1, \dots, x_{n+p}$, the orders of zeroes or poles 
being prescribed by the list
$k_1, \dots, k_{n+p}$ ($k_i >0$ for the zeroes, $k_i <0$ for the
poles, and $k_i=0$ for the marked points that are neither zeroes nor
poles). The closure of this set in $\oM_{g,n}$
is called the {\em double ramification
cycle} or a {\em DR-cycle}.
\end{definition}

Here are some basic facts about the double ramification cycles that
make them so useful.

\begin{enumerate}
\item
The codimension of a DR-cycle is equal to $g-p$ whenever there is at
least one positive and one negative number among $k_1, \dots, k_n$
(see~\cite{Mumford, Ionel}).
Assuming that this condition is satisfied we see that
for $p=g$ the DR-cycle coincides with the moduli space $\oM_{g,n}$.

\item
The cohomology class Poincar{\'e} dual to
any DR-cycle belongs to the tautological ring of~$\oM_{g,n}$ 
(proved in~\cite{FabPan}). 

This makes the results of this paper applicable to DR-cycles.

\item
Each DR-cycle is the image of the corresponding DR-space under
the forgetful map~$\pi$ that forgets the covering and all the
marked points except $x_1, \dots, x_n$, but retains and
stabilizes the source curve with the $n$ remaining marked points. 
The map~$\pi$ sends the fundamental
homology class of the DR-space to a multiple of the fundamental
homology class of the DR-cycle.

\item
Every class $\pi^*(\psi_i)$ on a DR-space can be expressed
as a linear combination of boundary divisors~\cite{Ionel,Shadrin1}.

This is a very important property that can be used to compute
integrals involving $\psi$-classes.  Indeed, it allows us to
get rid of the $\psi$-classes one by one by reducing the integral
to simpler integrals over smaller spaces. Using this procedure,
the following result was established in~\cite{ShaZvo}, Theorem~1.
\end{enumerate}

\begin{theorem} \label{Thm:rspin} {\rm \cite{ShaZvo}}
Every correlator in the $r$-spin Witten conjecture 
can be expressed in genus zero correlators using only 
tautological relations.
\end{theorem}

The last result shows that Theorem~\ref{Thm:GeomForm}
is applicable to the $r$-spin Witten conjecture and
suffices to prove it.

\bigskip

The proof of Theorem~\ref{Thm:rspin} in~\cite{ShaZvo}
goes as follows. As explained in Item~1, the fundamental
class of~$\oM_{g,n}$ can be represented as a DR-cycle
with $p=g$. If this is done in an intelligent way,
then the elimination of $\psi$-classes 
according to Item~4, leads us to boundary divisors that
can themselves be expressed as DR-cycles on the boundary.
We end up with the integral of Witten's class $c_W$ (without $\psi$-classes)
over a DR-cycle. A dimension count shows that an integral like that
may be nonzero only if $g=0$ or~1 and the codimension of
the DR-cycle is equal to the genus (0 for $g=0$ and 1 for $g=1$).
In the latter case
we must do some more work: putting one $\psi$-class back into
the integral and expressing it as a linear combination of 
boundary divisors in different ways we obtain certain relations
between genus~1 and genus~0 integrals of~$c_W$. It turns out that 
these relations suffice to reduce all genus~1 integrals to genus~0
integrals.

However, now we can give a simpler, although less constructive,
proof.

\paragraph{Proof of Theorem~\ref{Thm:rspin}.}

In~\cite{Ionel}, Ionel proved the following assertion:

{\em Let $M$ be a monomial in $\psi$- and $\kappa$-classes
on $\oM_{g,n}$, of degree\footnote{The degree of a
$\psi$-class equals~$1$, while the degree of $\kappa_k$ equals~$k$.} 
at least $g$ for $g \geq 1$ or at least~$1$ for $g=0$. 
Then the class~$M$ can be represented as a linear combination
of classes of the form
$$
q_* 
\left[
(\mbox{\rm DR-cycle})_1 \times \dots \times (\mbox{\rm DR-cycle})_k
\right].
$$
}
Here $k \geq 1$ is an integer that can be different 
for different terms of the sum, 
$q: \oM_{g_1,n_1} \times \dots \times \oM_{g_k,n_k}
\rightarrow \oM_{g,n}$ is the gluing map from a product of smaller
moduli spaces to a boundary stratum of $\oM_{g,n}$
and the cycles $(\mbox{\rm DR-cycle})_j$ are DR-cycles on
the smaller moduli spaces. 

It was established in~\cite{FabPan} that every DR-cycle is
tautological and that Ionel's theorem can therefore be
improved in the following way:

{\em Let $M$ be a monomial in $\psi$- and $\kappa$-classes
on $\oM_{g,n}$, of degree
at least $g$ for $g \geq 1$ or at least~$1$ for $g=0$. 
Then the class~$M$ can be represented by a linear combination
of dual graphs each of which has at least one edge.}

We will call this property the {\em $g$-reduction}.

Now, a simple dimension count shows that the integral
$$
\int\limits_{\oM_{g,n}} \!\!\!
\beta \cdot  c_W(a_1, \dots, a_n)
$$
vanishes unless the class $\beta$ has complex degree at least~$g$.
Indeed, the degree of Witten's class equals
$$
\deg c_W = \frac{(r-2)(g-1)+\sum a_i}r \leq \frac{(r-2)(n+g-1)}r,
$$
while the dimension of $\oM_{g,n}$ is $3g-3+n$
(use the exact expression for $g=1$ and the upper bound for
$g \geq 2$).

Recall that the pull-back of Witten's class to the boundary
components is given by the factorization property~(3)
of Section~\ref{Ssec:WiGrWi}. This makes it easy to apply
the $g$-reduction to integrals involving Witten's class.

The rest of the proof is simple. 
Suppose we wish to compute the integral
$$
\int\limits_{\oM_{g,n}} \!\!\!
c_W(a_1, \dots, a_n) \, \psi_1^{d_1} \dots \psi_n^{d_n}.
$$
Apply the $g$-reduction in iteration as many times as possible,
starting with the class $\psi_1^{d_1} \dots \psi_n^{d_n}$.
In the end we will obtain an expression of
$\prod_{i=1}^n \psi_i^{d_i}$ as a linear combination 
of dual graphs~$G$ satisfying the following condition.
Suppose a vertex~$v$ of~$G$ is labeled with genus $g_v >0$ and with a class
$\kappa_{k_1, \dots, k_m}$, and suppose the half-edges and tails
issuing from~$v$ are labeled with $\psi^{d_1}$,
\dots, $\psi^{d_{n_v}}$. Then $\sum k_i + \sum d_i < g_v$.
(Indeed, if $\sum k_i + \sum d_i \geq g_v$ for at least one vertex,
we can apply the $g$-reduction to this vertex.) 
But, as we have already explained,
the integral of Witten's class over the class $[G]$ represented 
by a dual graph like that vanishes whenever there is at least one 
vertex of nonzero genus. Thus the only contribution comes
from graphs with only genus zero vertices with no $\psi$- or
$\kappa$-classes. So we have
reduced any given correlator involved in Witten's conjecture
to a linear combination of products of genus~0 correlators.
More precisely, the only remaining correlators are
integrals of Witten's class with no $\psi$-classes
over genus zero moduli spaces.
\qed

\bigskip

As we have already explained, this implies that the formal
Gromov--Witten potential $F^{(r)}$ coincides with the
geometric Gromov--Witten potential  $F^{[r]}$ and proves
the Witten conjecture.

\section{More on Givental's quantization and Y.-P. Lee's theorem}
\label{Sec:Lee}

A proof of Theorem~\ref{Thm:Lee} is contained in~\cite{Lee2}.
However there are some missing details that we would like to
fill in here. First, we would like to explain
precisely why the operators $\tau_k$ act on the $\kappa$-classes
in the way described in Section~\ref{Sssec:tau}
(this is done in Proposition~\ref{Prop:Rderiv}).
Second, we explain more precisely how Givental's quantization
is applied to prove the theorem.

``Givental's group'' is really not a group, but a collection
of two groups: the so-called 
``lower triangular'' and ``upper triangular'' groups.
This is analogous to the Birkhoff decomposition in
the finite-dimensional case. However, because both
groups are infinite-dimensional, it is in general not
possible to multiply their elements, like it is impossible
to multiply a power series in~$z$ and a power series in $z^{-1}$.
On the other hand, it turns out that, under some conditions,
one can apply first an element of the upper triangular
group and then an element of the lower triangular group
to a Gromov--Witten potential.

Both the lower triangular and the upper triangular group
possess Lie algebras. An element of the lower triangular Lie algebra is a
series $s(z^{-1}) = \sum_{l \geq 1} s_l z^{-l}$ 
of linear operators on the vector
space~$V$. The operators $s_l$ are self-adjoint
for $l$ odd and skew-self-adjoint for $l$ even with respect
to the quadratic form~$\eta$. 

Similarly, an element of the upper triangular Lie algebra is a
series $r(z) = \sum_{l \geq 1} r_l z^l$ 
of linear operators on~$V$. The operators $r_l$ are self-adjoint
for $l$ odd and skew-self-adjoint for $l$ even with respect
to the quadratic form~$\eta$.

The coefficients of $s_l$ and $r_l$ in the basis~$A$ will be denoted by
$(s_l)_\mu^\nu$ and $(r_l)_\mu^\nu$. We will also need
the bivectors and the bilinear forms given by
$$
(s_l)_{\mu\nu} = \sum_\rho \eta_{\mu\rho} (s_l)^\rho_\nu,
\qquad
(s_l)^{\mu\nu} = \sum_\rho \eta^{\mu\rho} (s_l)^\nu_\rho,
$$
$$
(r_l)_{\mu\nu} = \sum_\rho \eta_{\mu\rho} (r_l)^\rho_\nu,
\qquad
(r_l)^{\mu\nu} = \sum_\rho \eta^{\mu\rho} (r_l)_\rho^\nu.
$$
The matrices $(s_l)_{\mu\nu}$, $(r_l)_{\mu\nu}$, 
$(s_l)^{\mu\nu}$, and  $(r_l)^{\mu\nu}$ are
symmetric for odd~$l$ and skew-symmetric for even~$l$.

Y.-P.~Lee~\cite{Lee2} writes down explicit formulas for the action
of $s$ and $r$ on any given correlator of a Gromov--Witten
potential~$F$ (see below). Once this is done, the main problem is to
understand what happens when we apply these formulas to the
induced vector of a tautological class: indeed, both
Y.-P.~Lee's formulas for the derivatives of an individual
correlator and the expression of the induced vector
in terms of correlators (described in
Section~\ref{Ssec:univrel}) are fairly complicated.

Below we sum up the argument of~\cite{Lee2}
and give more detailed statements of certain results.

If $s(z^{-1}) = \sum_{l \geq 1} s_l z^{-l}$ is an element
of the lower triangular Lie algebra and 
$r(z) = \sum_{l \geq 1} r_l z^l$ an element of the upper
triangular Lie algebra, denote, for shortness
$$
s_l(\tau_{d,\mu}) = \sum_\nu (s_l)^\nu_\mu \; \tau_{d-l,\nu},
\qquad
r_l(\tau_{d,\mu}) = \sum_\nu (r_l)^\nu_\mu \; \tau_{d+l,\nu}.
$$

Now we are going to follow the path from the Gromov--Witten
potential of a point to the general formal semi-simple Gromov--Witten
potential.

\subsection{The upper triangular group}
\label{Ssec:uptriang}

Let $t_0^\mu$ be flat coordinates on a semi-simple 
Frobenius manifold $M$ of dimension~$k$. Let $f(t_0^\mu)$ be the 
corresponding genus~0 potential. We assume that
$f$ is an analytic function. At a semi-simple point,
the tangent Frobenius algebra~$T_*M$ to
the Frobenius manifold~$M$ possesses a basis of primitive
idempotents. Denote by $\alpha_1, \dots, \alpha_k$ their 
scalar squares in the metric~$\eta$. Givental constructs an element~$R$
of the upper triangular group whose action transforms
the constant Frobenius structure on $T_*M$ into the Frobenius
structure of~$M$ at the neighborhood of the semi-simple 
point.

In other words, the first step of Givental's quantization
is to apply the quantized action of $R$ to
the partition function
$Z^{\pt, \alpha_1} \times \dots \times Z^{\pt, \alpha_k}$,
see Eq.~(\ref{Eq:rescaled}).

The potential $F^{\pt, \alpha_1} \oplus \dots \oplus F^{\pt, \alpha_k}$ 
obviously possesses the two following crucial properties.

\begin{definition}
A Gromov--Witten potential is {\em tame} if
$\left< \tau_{d_1}^{\mu_1} \cdots \tau_{d_n}^{\mu_n}\right>_g $
vanishes whenever $\sum d_i > 3g-3+n$.
\end{definition}

\begin{definition}
A Gromov--Witten potential is an {\em ancestor} potential
if its correlators with $2-2g-n \geq 0$ vanish.
\end{definition}

Let $r = \ln R$ be an element of the upper triangular Lie algebra.
In general, the action of $R$ on a power series is not well defined.
However, it is easy to check (cf.~Proposition~\ref{Prop:Rformula})
that the action of $r_l$ increases the grading $\sum d_i -n -3g+3$
by~$l$. Therefore the action of $R$ on a tame series is
well-defined and is equal to the exponential of the action of~$r$.

According to Givental's formulas, $r = \sum_{l \geq 1} r_l z^l$ 
acts on the partition function~$Z$ via the second 
order differential operator
$$
\widehat{r} = \hspace{20cm}
$$
$$
- \sum_{
\substack{l \geq 1 \\ \mu}
}
(r_l)_1^\mu \frac{\d}{\d t_{l+1}^\mu}
+\sum_{
\substack{d \geq 0, l \geq 1\\ \mu,\nu}
} \!\! 
(r_l)_\nu^\mu \, t_d^\nu \frac{\d}{\d t_{d+l}^\mu}
+\frac{\hbar}2 \sum_{
\substack{d_1,d_2 \geq 0 \\ \mu_1,\mu_2}
} \!
(-1)^{d_1+1} (r_{d_1+d_2+1})^{\mu_1 \mu_2}
\frac{\d^2}{\d t_{d_1}^{\mu_1} \d t_{d_2}^{\mu_2}}.
$$

The next proposition gives the action of~$r$ on individual correlators.

\begin{proposition} {\rm (\cite{Lee2}, Equation~(23))}
\label{Prop:Rformula}
\begin{eqnarray}
\label{Eq:Rformula}
\lefteqn{r. \left< \tau_{d_1, \mu_1} \dots \tau_{d_n, \mu_n} \right>_g =} \\
\nonumber
&&{}- \sum_{l=1}^{\infty} 
\left<
r_l(\tau_{1,1}) \, \tau_{d_1, \mu_1} \dots \tau_{d_n,\mu_n}
\right>_g\\
\nonumber
&&
{}+ \sum_{l=1}^{\infty} \sum_{i=1}^n
\left<
\tau_{d_1, \mu_1} \dots r_l(\tau_{d_i,\mu_i}) \dots \tau_{d_n,\mu_n}
\right>_g\\
\nonumber
&&{}+\frac12 \sum_{l=1}^{\infty} \sum_{m+m' = l-1} \!\!\!\! (-1)^{m+1}
\sum_{\mu,\nu} (r_l)^{\mu\nu}
\left<
\tau_{m,\mu} \tau_{m',\nu} \, \tau_{d_1,\mu_1} \dots \tau_{d_n,\mu_n}
\right>_{g-1}\\
\nonumber
&&{}+\frac12 \sum_{l=1}^{\infty} \sum_{m+m'=l-1} \!\!\!\! (-1)^{m+1}
\!\!\!\!\!\!\!\!\!
\sum_{\stackrel{\scriptstyle g_1+g_2 = g}
{I \sqcup J = \{1, \dots, n \}}
} \!\!\!
\sum_{\mu,\nu}
(r_l)^{\mu\nu}
\,
\biggl< \!
\tau_{m,\mu} \prod_{i \in I} \tau_{d_i,\mu_i}
\! \biggr>_{\!\!g_1}
\!
\biggl<\!
\tau_{m',\nu} \prod_{i \in J} \tau_{d_i,\mu_i}
\!\biggr>_{\!\!g_2} \!\!\!.
\end{eqnarray}
\end{proposition}

\begin{remark}
One can easily deduce from these formulas that
the action of the upper triangular group preserves
the tameness property and the property of being an
ancestor potential.
\end{remark}

\begin{proposition} {\rm (based on \cite{Lee2}, Section~6)}
\label{Prop:Rderiv}
Let $L$ be a linear combination of dual graphs and $r$ an element 
of the upper triangular Lie subalgebra of Givental's Lie algebra.
Let $F$ be a formal tame ancestor potential. Then
$$
r . \biggl<
p^*(L) \prod_{i=1}^{n+n'} \psi_i^{d_i} \ev_i^*(\mu_i)
\biggr>_g
$$
is a linear combination of polynomials in correlators of
the form
$$
{\rm (i)} \quad 
\biggl<
p^*(L) \prod_{i=1}^{n+n'+1} \psi_i^{d_i} \ev_i^*(\mu_i)
\biggr>_g ,
$$
where $d_{n+n'+1} \geq 2$ and $\mu_{n+n'+1}$ is a primary field;
$$
{\rm (ii)} \quad 
\biggl<
p^*(L) \prod_{i=1}^{n+n'} \psi_i^{d'_i} \ev_i^*(\mu'_i)
\biggr>_g ,
$$
where $d'_i = d_i$, $\mu'_i = \mu_i$ for all $i$ except one,
and $d'_i > d_i$ for one~$i$;
$$
{\rm (iii)} \quad 
\biggl<
p^*(\tau_l(L)) \prod_{i=1}^{n+n'} \psi_i^{d_i} \ev_i^*(\mu_i) \cdot
\ev_\alpha^*(\mu_\alpha) \ev_\beta^*(\mu_\beta) 
\biggr>_g ,
$$
where $l \geq 1$  and $\mu_\alpha, \mu_\beta$ are primary fields;
$$
{\rm (iv)} \quad 
\biggl<\tau_{d_\alpha,\mu_\alpha} \prod_{i \in I} \tau_{d_i, \mu_i}
\biggr>_{\!\!0}
\biggl<p^*(L) \; \psi_\beta^{d_\beta} \ev_\beta^*(\mu_\beta)
\prod_{i \in J} \psi_i^{d_i} \ev_i^*(\mu_i)
\biggr>_{\!\! g} ,
$$
where $I \subset \{ n+1, \dots, n+n' \}$,
$I \sqcup J = \{ 1, \dots, n+n' \}$, $d_\alpha$ and $d_\beta$
are nonnegative integers,
and $\mu_\alpha, \mu_\beta$ are primary fields;
$$
{\rm (v)} \quad 
\biggl<\tau_{d_\alpha,\mu_\alpha} \tau_{d_j,\mu_j} 
\prod_{i \in I} \tau_{d_i, \mu_i} \biggr>_{\!\! 0}
\biggl<
p^*(L) \, \psi_\beta^{d_\beta} \ev_\beta^*(\mu_\beta)
\prod_{i \in J} \psi_i^{d_i} \ev_i^*(\mu_i)
\biggr>_{\!\! g} ,
$$
where $I \subset \{ n+1, \dots, n+n' \}$, $j \in \{ 1, \dots, n \}$,
$I \sqcup J \sqcup \{ j \} = \{ 1, \dots, n+n' \}$,
$d_\alpha$ and $d_\beta$ are nonnegative integers,
and $\mu_\alpha, \mu_\beta$ are primary fields.
\end{proposition}

\begin{remark}
If we work over a Novikov ring, we must introduce a
summation over $D_1 + D_2 = D$ in the last term of the
equality in Proposition~\ref{Prop:Rformula} as well as in
terms~(iv) and~(v) of Proposition~\ref{Prop:Rderiv}.
\end{remark}

\begin{corollary} \label{Cor:Rinv}
Let $L$ be a linear combination of genus~$g$ dual graphs with
$n$ tails such that
$\tau_k(L)$ is a universal relation for all $k$. Assume
that $F^\pt_L$ vanishes. Then $F_L$ also vanishes for every
$F$ that can be obtained from 
$F^{\pt, \alpha_1}_0 \oplus \dots \oplus F^{\pt, \alpha_k}_0$ 
by the upper triangular group action.
\end{corollary}

\paragraph{Proof.} The polynomials (i) to~(v) are either
entries of~$F_L$ of~$F_{\tau_l(L)}$. In the latter case they vanish
by assumption. Therefore the vector $F_L$ on the orbit
$e^{tr} F$ is the solution of a linear
differential equation with vanishing initial conditions. Thus
$F_L=0$ on the whole orbit. \qed

\paragraph{Proof of Proposition~\ref{Prop:Rderiv}.} \ \\
{\bf 1.} We first prove the proposition in the
particular case $n'=0$, $d_1 = \dots = d_n=0$. 
In Section~\ref{Ssec:univrel} we gave a three steps algorithm 
to determine the linear combination of graphs
$p^*(L) \prod_{i=1}^{n+n'} \psi_i^{d_i} \ev_i^*(\mu_i)$.
In our particular case, the steps are greatly simplified:
Step~1 (replacing $\kappa$-classes by additional tails)
remains unchanged; Step~2 (expressing ancestor $\psi$-classes
in terms of descendant $\psi$-classes) can be skipped, since
all genus zero 2-point correlators vanish for an ancestor
potential and since there are no additional marked points; Step~3
(adding tails from $n+1$ to $n+n'$ and multiplying by
$\psi_i^{d_i} \ev_i^*(\mu_i)$) reduces to assigning
the markings $\mu_1, \dots, \mu_n$ to tails $1$ to~$n$.

As for the claim of our proposition, if $n'=0$, then 
terms~(iv) and~(v) disappear, because 
$I \subset \{ n+1, \dots, n+n' \}$ is then empty, and an ancestor
potential does not have genus~0 correlators with fewer than~3
entries. Thus we must prove that 
$r . \left<
p^*(L) \prod_{i=1}^{n} \ev_i^*(\mu_i)
\right>_g$
is a linear combination of terms~(i), (ii), and~(iii).

Let~$G$ be a graph in~$L$.

{\bf 1a.} First suppose that $G$ has only one vertex (of genus~$g$
and valency~$n$) with a class $\kappa_{k_1, \dots, k_m}$
assigned to the vertex and classes $\psi_i^{d_i}$ assigned to
the tails. Applying Steps~1 and~3 we obtain 
$$
\left<
p^*(G) \prod_{i=1}^{n} \ev_i^*(\mu_i)
\right> = \biggl< \prod_{i=1}^m \tau_{k_i+1,1} \;
\prod_{i=1}^{n} \tau_{d_i,\mu_i} \biggr>.
$$
To determine the action of $r$ on the correlator in the
right-hand side we apply Equation~\eqref{Eq:Rformula} to it.

Applying the third and fourth terms in Equation~(\ref{Eq:Rformula})
gives us term~(iii) involving $\tau_l(G)$. 
Note that when we apply the fourth term 
of Equation~(\ref{Eq:Rformula}), the indices $k_1, \dots, k_m$
are distributed among the two correlators in all
possible ways, according to the description of $\tau_l$.

Now apply the first two terms of Equation~(\ref{Eq:Rformula})
to our correlator. We regroup them in the following way:
$$
\sum_{l=1}^{\infty} \sum_{\mu} (r_l)^1_\mu
\left[\sum_{j=1}^m \biggl< \tau_{k_j+l+1,\mu} \prod_{i \not= j} 
\tau_{k_i+1,1} \prod_{i=1}^{n} \tau_{d_i,\mu_i} \biggr>
-
\biggl< \tau_{l+1,\mu} \prod_{i=1}^m \tau_{k_i+1,1}
\prod_{i=1}^{n} \tau_{d_i,\mu_i} \biggr>\right] 
$$
\begin{equation} \label{Eq:wellbeing}
{} +
\sum_{j=1}^{n} \biggl< r_l(\tau_{d_j,\mu_j}) \prod_{i=1}^m \tau_{k_i+1,1}
\prod_{i \not= j} \tau_{d_i,\mu_i} \biggr>.
\end{equation}

Consider the forgetful map $\pi:\oM_{g,n+1} \to \oM_{g,n}$.
The pull-back $\pi^*(\kappa_{k_1, \dots, k_m})$ 
given by Lemma~\ref{Lem:kappa} imitates the expression in
square brackets of~\eqref{Eq:wellbeing}. 
In addition, we have $\pi^*(\psi_i^{d_i}) \psi_{n+1}^{l+1} = 
\psi_i^{d_i} \psi_{n+1}^{l+1}$. It follows that 
the terms in square brackets add up to give term~(i).

The last term of Equation~(\ref{Eq:wellbeing})
replaces one of the symbols $\tau_{d_j,\mu_j}$
by $\tau_{d_j',\mu_j'}$ with $d'_j > d_j$, which
gives us term~(ii).

{\bf 1b.} Now let $G$ be an arbitrary dual graph. Then
the terms in square brackets of 
Equation~(\ref{Eq:wellbeing}) will appear for
each vertex of the graph and their sum will still represent the
class $\pi^*(G) \cdot\psi_{n+1}^{l+1} \ev_{n+1}^*(\mu)$. 

The last term of Equation~(\ref{Eq:wellbeing})
replaces as before, one of the symbols
$\tau_{d_j,\mu_j}$ by $\tau_{d_j',\mu_j'}$. But now this symbol
can either correspond to a tail or to a half-edge of~$G$. 
If the symbol we replace corresponds to a tail, 
it gives rise to term~(ii) in the
proposition. If it corresponds to a half-edge of~$G$, 
the corresponding term 
$\biggl< r_l(\tau_{d_j,\mu_j}) \prod_{i=1}^m \tau_{k_i+1,1}
\prod_{i \not= j} \tau_{d_i,\mu_i} \biggr>$
contributes to term~(iii) of the proposition. Finally, the third 
and fourth terms in Equation~(\ref{Eq:Rformula})
complete the expression for term~(iii).

{\bf 2.} Now we return to the general case, when $n'$
and $d_1, \dots, d_{n+n'}$ are arbitrary. Consider the
forgetful maps
$$
\oM_{g,n+n'+1}
\stackrel{\pi_1}{\rightarrow}
\oM_{g,n+n'}
\stackrel{\pi_2}{\rightarrow}
\oM_{g,n}
$$
and the composition $\pi = \pi_2 \circ \pi_1$.

If~$L$ is our initial linear
combination of dual graphs, then
$L' = \pi_2^*(L) \prod_{i=1}^{n+n'} \psi_i^{d_i}$ is a well-defined
linear combination of dual graphs obtained by the usual
three steps algorithm, but without attaching the markings $\mu_i$
to the tails.

Now we can apply the particular case that we have just proved
to the linear combination~$L'$. All we have to do is 
reinterpret the answer in terms of~$L$. We claim the following.

Term~(i) for~$L'$ gives term~(i) for~$L$. Indeed, the 
classes on the moduli spaces of curves involved in this term
are 
$$
\pi_1^*(L') \;  \psi_{n+n'+1}^{d_{n+n'+1}} \;=\; 
\pi^*(L) \;
\pi_1^*\biggl(\prod_{i=1}^{n+n'} \psi_i^{d_i} \biggr) \;
\psi_{n+n'+1}^{d_{n+n'+1}}.
$$
But since $d_{n+n'+1} \geq 2$, we have
$\pi_1^*(\psi_i^{d_i}) \psi_{n+n'+1}^{d_{n+n'+1}} = 
\psi_i^{d_i} \psi_{n+n'+1}^{d_{n+n'+1}}$ on $\oM_{g,n+n'+1}$.

Term~(ii) for~$L'$ gives term~(ii) for~$L$. This is obvious.

Term~(iii) for~$L'$ gives terms~(iii), (iv), and~(v) for~$L$.
Indeed, the space $\oMbd_{g,n+n'}$ has two kinds of irreducible
components: those that project to a component of
$\oMbd_{g,n}$ under $\pi_2$ and those that project
onto $\oM_{g,n}$ under~$\pi_2$. In the first case
we obtain term~(iii) of the proposition.
As to the second case, it occurs
when the generic curve of the component of $\oMbd_{g,n+n'}$
has one component of genus~$g$ and one component of genus~$0$
that is contracted by~$\pi_2$. The class $\pi_2^*(L)$ is supported
on the genus~$g$ component. The contracted component of genus~$0$
contains either~0 or~1 point with markings $1,\dots,n$
(those that are not forgotten by~$\pi_2$). According to these
two cases we obtain either term~(iv) or term~(v) from
the proposition.
\qed

\subsection{The lower triangular group}

By the action of the upper triangular group we
have obtained a Gromov--Witten potential that describes
the semi-simple Frobenius manifold~$M$ that we started with.
Recall that this Gromov--Witten potential is
a power series in variables $t_1^\mu, t_2^\mu, \dots$
whose coefficients are functions in variables
$t_0^\mu$ analytic outside the discriminant of~$M$.

Now the action of the lower triangular group 
re-expands the Gromov--Witten potential at a
different (possibly non semi-simple) point of~$M$ and simultaneously 
incorporates 1- and 2-point genus~0 correlators.

The element $s = \sum_{l \geq 1} s_l z^{-l}$ acts on the partition
function~$Z$ via the first order differential operator
$$
\whs = - \sum_\mu (s_1)_1^\mu \frac{\d}{\d t_0^\mu} \hspace{20cm}
$$
$$
+ \frac1\hbar \sum_{d,\mu} (s_{d+2})_{1,\mu} \, t_d^\mu
+ \sum_{
\substack{d,l\\ \mu, \nu}
}
(s_l)_\nu^\mu \, t_{d+l}^\nu \frac{\d}{\d t_d^\mu}
+ \frac1{2 \hbar} \sum_{
\substack{d_1,d_2 \\ \mu_2,\mu_2}
}
(-1)^{d_1} (s_{d_1+d_2+1})_{\mu_1,\mu_2} \, t_{d_1}^{\mu_1} t_{d_2}^{\mu_2}.
$$

In this expression we have omitted the term $-\frac1{2\hbar} (s_3)_{1,1}$.
Indeed, it commutes with all other terms and the action of
its exponential only adds the constant $(s_3)_{1,1}$ to $F_0$. 
Similarly, in the
sequel we will consider $F_0$ to be defined up to an additive 
constant and omit those terms in differential operators that
do nothing more than changing this constant.

The following proposition makes the action of $s$ on
individual correlators explicit.

\begin{proposition} {\rm (\cite{Lee2}, Equations~(19),(20))}
\label{Prop:Sformula}
For $2-2g-n <0$, we have
\begin{eqnarray}
\label{Eq:Sformula1}
\lefteqn{s.\left< \tau_{d_1,\mu_1}\dots\tau_{d_n,\mu_n} \right>_g=} \\
\nonumber
&&- \left<s_1(\tau_{1,1})\, \tau_{d_1,\mu_1} \dots \tau_{d_n,\mu_n} \right>_g
+ \sum_{l=1}^{\infty} \sum_{i=1}^n 
\left<
\tau_{d_1,\mu_1} \dots s_l(\tau_{d_i,\mu_i}) \dots\tau_{d_n,\mu_n}
\right>_g.
\end{eqnarray}
If $g = 0$ and $n = 2$, we have
\begin{eqnarray}
\label{Eq:Sformula2}
\lefteqn{s . \left< \tau_{d_1,\mu_1} \tau_{d_2,\mu_2} \right>_0 =} \\
\nonumber
&&-\left<s_1(\tau_{1,1})\, \tau_{d_1,\mu_1} \tau_{d_2,\mu_2} \right>_0
+ \sum_{l=1}^{\infty}
\left[
\left<s_l(\tau_{d_1,\mu_1}) \tau_{d_2,\mu_2}\right>_0
+\left<\tau_{d_1,\mu_1} s_l(\tau_{d_2,\mu_2})\right>_0
\right]\\
\nonumber
&&{}+(-1)^{d_1} (s_{d_1+d_2+1})_{\mu_1\mu_2}.
\end{eqnarray}
If $g=0$ and $n=1$, we have
\begin{eqnarray}
\label{Eq:Sformula3}
s . \left< \tau_{d,\mu}\right>_0 &=&
-\left<s_1(\tau_{1,1})\, \tau_{d,\mu} \right>_0
+ \sum_{l=1}^{\infty}
\left<s_l(\tau_{d,\mu})\right>_0
+ (s_{d+2})_{1\mu}.
\end{eqnarray}
\end{proposition}

Looking at Equations~(\ref{Eq:Sformula1}), 
(\ref{Eq:Sformula2}), (\ref{Eq:Sformula3}) 
we see that $s_l$ decreases the grading $\sum d_i$ by~$l$,
except for the first term involving $s_1(\tau_{1,1})$.
This annoying term leads to a problem that we have 
to discuss in more detail. The corresponding term
$$
\whu= \sum_\mu (s_1)^\mu_1 \frac{\d}{\d t_0^\mu}
$$ 
in the differential operator $\widehat s$ is
simply a partial derivative in the direction $s_1(1)$.
The exponential $\exp(\whu)$ is then a
shift of the coordinates $t_0^\mu$ in this direction.
A shift of coordinates is not a well-defined
operation for formal power series. In our case, we
consider power series in variables $t_1^\mu$, $t_2^\mu$, \dots,
whose coefficients are analytic in $t_0^\mu$ on the
Frobenius manifold except perhaps its discriminant.
Thus the shift is a well-defined operation, but
it can, in some cases, lead us out of the realm of
power series. In this case the genus expansion of
the genus~0 Gromov--Witten potential will be a
power series in $t_1^\mu$, $t_2^\mu$, \dots, 
whose coefficients are analytic functions in $t_0^\mu$
with a singularity at the origin. 
(See, for instance,~\cite{DubZha}, Section~6, where
$F$ and $G$ are the Gromov--Witten potentials without
descendants for genus $0$ and~$1$ respectively.)

We also need to give a proper definition of the action
of $S = \exp(s)$. We have, up to omitted constant terms, 
\begin{eqnarray*}
[\whs, \whu] &=&
\frac1\hbar \sum_{d,\mu} (s_1 s_{d+1})_{1,\mu} \,  t_d^\mu,\\
\left[\whs,[\whs, \whu]\right] &=&
\frac1\hbar \sum_{l,d,\mu} (s_l s_1 s_{d+1})_{1,\mu} \,  t_{d+l}^\mu,\\
\left[\whs,[\whs,[\whs, \whu]]\right] &=&
\frac1\hbar \sum_{l_1, l_2,d,\mu} 
(s_{l_1} s_{l_2} s_1 s_{d+1})_{1,\mu} \,  t_{d+l_1+l_2}^\mu,
\end{eqnarray*}
etc. On the other hand, $[\whu,[\whu, \whs]]=0$
(again, up to constant term). Thus, according to the 
Baker--Campbell--Hausdorff formula,
we have $e^{\whs} = e^{-\whu} e^{\whv}$, where
\begin{equation} \label{Eq:v}
\whv = \whu+\whs - \frac12 [\whs, \whu] + \frac1{12} 
[\whs,[\whs, \whu]] - \frac1{720} [\whs,[\whs,[\whs, \whu]]] + \dots
\end{equation}
$$
= \sum_{
\substack{d,l\\ \mu, \nu}
}
(s_l)_\nu^\mu \, t_{d+l}^\nu \frac{\d}{\d t_d^\mu}
+ \frac1{2 \hbar} \sum_{
\substack{d_1,d_2 \\ \mu_2,\mu_2}
}
(-1)^{d_1} (s_{d_1+d_2+1})_{\mu_1,\mu_2} \, t_{d_1}^{\mu_1} t_{d_2}^{\mu_2}
+ \sum C_{d,\mu} t_d^\mu,
$$
where every $C_{d, \mu}$ is a finite polynomial in the matrix elements
$(s_l)_{\mu\nu}$. 
\begin{proposition}
\label{Prop:Uformula}
For $2-2g-n <0$, we have
\begin{equation}
\label{Eq:Uformula1}
\whv.\left< \tau_{d_1,\mu_1}\dots\tau_{d_n,\mu_n} \right>_g =
\sum_{l=1}^{\infty} \sum_{i=1}^n 
\left<
\tau_{d_1,\mu_1} \dots s_l(\tau_{d_i,\mu_i}) \dots\tau_{d_n,\mu_n}
\right>_g.
\end{equation}
If $g = 0$ and $n = 2$, we have
\begin{eqnarray}
\label{Eq:Uformula2}
\lefteqn{\whv . \left< \tau_{d_1,\mu_1} \tau_{d_2,\mu_2} \right>_0 =} \\
\nonumber &&
\sum_{l=1}^{\infty}
\left[
\left<s_l(\tau_{d_1,\mu_1}) \tau_{d_2,\mu_2}\right>_0
+\left<\tau_{d_1,\mu_1} s_l(\tau_{d_2,\mu_2})\right>_0
\right]+(-1)^{d_1} (s_{d_1+d_2+1})_{\mu_1\mu_2}.
\end{eqnarray}
\end{proposition}

\begin{remark}
The action of $\whv$ on genus~0 one point correlators involves
the coefficients $C_d^\mu$ and is much more complicated.
However it is easy to see that such correlators never
appear in the induced vectors of tautological relations.
Therefore the action of $\whv$ on them is immaterial to us.
\end{remark}

\begin{lemma}
The action of $e^{\whv}$ is well-defined on power series;
in other words, every coefficient in the series $e^{\whv} Z$
is a finite polynomial in coefficients of~$Z$.
\end{lemma}

\paragraph{Proof.} The action of $\whv$ on a correlator is a
polynomial involving correlators with {\em strictly smaller} $\sum d_i$
and matrix elements of the matrices $s_l$. \qed

\begin{definition}
The action of $e^{\whs}$ on $Z$ is defined as the action
of $e^{\whv}$ followed by a translation of coordinates
$t_0^\mu$ by the vector $-(s_0)_1^\mu$.
\end{definition}

To sum up: the same element $s$ of the lower triangular Lie subalgebra
of Givental's Lie algebra determines {\em three} differential operators:
$\whs$, $\whu$, and $\whv$. They have the following properties:
the action of $e^{\whv}$ is well-defined on power series; the
action of $e^{\whu}$ is a translation of coordinates;
we have $e^{\whs} = e^{-\whu} \, e^{\whv}$.

%

\bigskip

Now let $L$ be a linear combination of dual graphs of genus~$g$
with $n$ tails and $\left<
p^*(L) \prod_{i=1}^{n+n'} \psi_i^{d_i} \ev_i^*(\mu_i)
\right>_g$ an entry of its induced vector. Let $s$
be an element of the lower triangular Lie subalgebra of
Givental's Lie algebra and $\whv$ the operator defined
by Equation~\eqref{Eq:v}.

\begin{proposition}
\label{Prop:Sderiv} The result of the action
$$
\whv.  \biggl<
p^*(L) \prod_{i=1}^{n+n'} \psi_i^{d_i} \ev_i^*(\mu_i)
\biggr>_g
$$
is a linear combination of polynomials in correlators of
the form
\begin{equation} \label{Eq:uL}
\biggl<
p^*(L) \prod_{i=1}^{n+n'} \psi_i^{d'_i} \ev_i^*(\mu'_i)
\biggr>_g
\end{equation}
with $d'_i = d_i$, $\mu'_i = \mu_i$ for all~$i$ except one,
while $d'_i < d_i$ for one~$i$.
\end{proposition}

\begin{corollary} \label{Cor:Sinv}
If the induced vector $F_L$ vanishes
for one formal Gromov--Witten potential~$F$ then it
also vanishes for all potentials that can be obtained
from~$F$ by the action of the lower triangular subgroup
of Givental's group.
\end{corollary}

\paragraph{Proof.} The polynomials~\eqref{Eq:uL} are themselves
entries of~$F_L$. Therefore the vector $F_L$ on the orbit
$e^{t\whv} F$ is the solution of a linear
differential equation with vanishing initial conditions. Thus
$F_L=0$ on the whole orbit. 

We know that the vanishing of $F_L$ 
can be expressed as a family of partial differential equations with constant coefficients.
These equations are preserved by translations of $t_0^\mu$.
It follows that the condition $F_L=0$ is preserved by the action
of the lower triangular group. \qed

\paragraph{Proof of Proposition~\ref{Prop:Sderiv}.}
We must compute
$$
\whv . \biggl<
p^*(L) \prod_{i=1}^{n+n'} \psi_i^{d_i} \ev_i^*(\mu_i)
\biggr>_g.
$$
The main element of the proof is the following observation:
every contribution of the second term in Equation~(\ref{Eq:Uformula2})
cancels with some contribution of the first term
in Equations~(\ref{Eq:Uformula1}) and~(\ref{Eq:Uformula2}).

Indeed, recall that in the definition of the induced vector $F_L$
every half-edge and every tail of each graph in~$L$ was replaced
by a ``stick'' of several edges (Step~2), and then new tails were
added in all possible ways to the new graphs (Step~3).
Consider the linear combination~$L_3$ of graphs obtained 
after Step~3 of the procedure.
In $L_3$ take two graphs $\Gamma_1$ and $\Gamma_2$
differing only in one fragment:
\begin{center}
\
\begin{picture}(200,90)
\put(0,0){\includegraphics{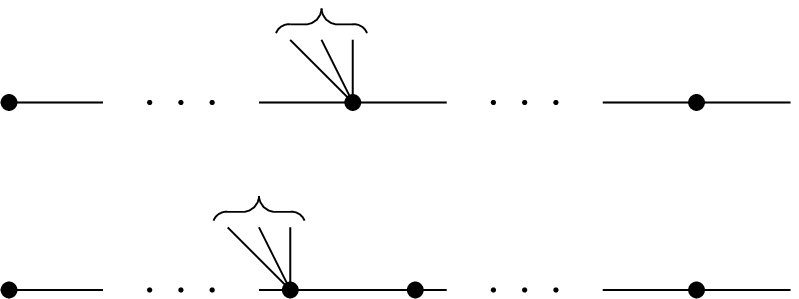}}
\put(107,61){$\psi^{d_a+d_b+1}$}
\put(125,7){$\psi^{d_a}$}
\put(90,7){$\psi^{d_b}$}
\put(-35,54){$\Gamma_1=$}
\put(-35,0){$\Gamma_2=$}
\put(70,32){$I$}
\put(88,86){$I$}
\end{picture}
\end{center}
Here $I \subset \{ n+1, \dots, n+n' \}$
is a set of labels of the tails added in Step~3.

The polynomial $P_{\Gamma_1}$ contains the factor
$$
Q_{\Gamma_1} =
\Bigl<\tau_{d_a+d_b+1,\mu_b} \tau_{0,\nu_b}
\prod_{i \in I} \tau_{d_i,\mu_i}
\Bigr>_{\! 0} \! .
$$ 

The polynomial $P_{\Gamma_2}$ contains the factor
$$
Q_{\Gamma_2} = 
\sum_{\mu_a,\nu_a} \left<\tau_{d_a,\mu_a} \tau_{0,\nu_a}\right>_0
\eta^{\nu_a\mu_b} \biggl<\tau_{d_b,\mu_b} \tau_{0,\nu_b}
\, \prod_{i \in I} \tau_{d_i,\mu_i}
\biggr>_0.
$$ 

Apply the operator~$\whv$ to the correlator
$\left<\tau_{d_a,\mu_a} \tau_{0,\nu_a}\right>_0$ in the
expression $Q_{\Gamma_2}$, and take the
contribution of the second term in Equation~(\ref{Eq:Uformula2}).
We obtain
$$
(-1)^{d_a} \sum_{\mu_a,\nu_a} (s_{d_a+1})_{\mu_a\nu_a} \;
\eta^{\nu_a\mu_b} \, \Bigl<\tau_{d_b,\mu_b} \tau_{0,\nu_b}
\prod_{i \in I} \tau_{d_i,\mu_i}
\Bigr>_{\! 0} \hspace*{10em}
$$
$$
\hspace*{10em}
= \Bigl<s_{d_a+1}(\tau_{d_a+d_b+1,\mu_b}) \, \tau_{0,\nu_b}
\prod_{i \in I} \tau_{d_i,\mu_i}
\Bigr>_{\! 0} \! .
$$
But this term is part of the action of $\whv$ on $P_{\Gamma_1}$. 
Indeed, it is part of the sum in the first term
of Equations~(\ref{Eq:Uformula1}) and~(\ref{Eq:Uformula2})
applied to
$$
\Bigl<\tau_{d_a+d_b+1,\mu_b} \tau_{0,\nu_b}
\prod_{i \in I} \tau_{d_i,\mu_i}
\Bigr>_{\! 0} \! .
$$
Since $\Gamma_2$ has one edge more than $\Gamma_1$, 
these terms appear with opposite
signs in Substitution~(1) and hence cancel.

Let us look at the contributions of the first term
of Equations~(\ref{Eq:Uformula1}) and~(\ref{Eq:Uformula2})
that survived the cancellation. They are exactly those
where $s_l$ is applied to the symbols $\tau_{d,\mu}$
corresponding to the $n+n'$ marked points and such that
$l$ is smaller than or equal to the corresponding~$d_i$.
These contributions combine into the expression given in the
proposition.
\qed

\paragraph{Proof of Theorem~\ref{Thm:Lee}.} The result
follows from Corollaries~\ref{Cor:Rinv} and~\ref{Cor:Sinv}.
\qed

\begin{remark}
Unfortunately, for the time being we cannot prove the implication
($L$ is a universal relation) $\Rightarrow$ ($\tau_1(L)$ is a universal
relation). There are two reasons for that. First, as explained in
Remark~\ref{Rem:UnivDisconnected}, for not necessarily connected
graphs checking that $Z_{\tau_1(L)}=0$ does not imply that $\tau_1(L)$
is a universal relation. Second, term~(iii) of 
Proposition~\ref{Prop:Rderiv} does not allow us to multiply
$\tau_1(L)$ by arbitrary powers of $\psi_\alpha$ and
$\psi_\beta$, but only by some of their combinations, namely,
the classes $\rho_k$ of Section~\ref{Ssec:tau}.
\end{remark}

\section*{Acknowledgements}
The authors are deeply grateful to M.~Kazarian for his
explanations on Frobenius manifolds and Givental's quantization.
We thank B.~Dubrovin for pointing out the necessity
of analyticity conditions on Gromov--Witten potentials.
The third author benefited greatly from a discussion with Tom Coates.
We also wish to thank the participants of the
Moduli Spaces program at the Mittag-Leffler Institute 
(Djursholm, Sweden) for the stimulating atmosphere
and lots of interesting discussions. We thank the
Institute for hospitality and support. 

C.F.~is supported by the grants 622-2003-1123 from the Swedish Research
Council and DMS-0600803 from the National Science Foundation
and by the G\"oran Gustafsson Foundation for Research in
Natural Sciences and Medicine.
D.Z.~is partly supported by the ANR project ``Geometry and
Integrability in Mathematical Physics'' ANR-05-BLAN-0029-01.

\end{document}